\documentclass[journal]{IEEEtran}
\IEEEoverridecommandlockouts                              

 \usepackage{etoolbox}
\makeatletter
\patchcmd{\@begintheorem}{\textit}{\textbf}{}{}
\makeatother
 \newtheorem{definition}{\bf Definition}
 \newtheorem{problem}{\bf Problem}

  \newtheorem{thm}{\bf Theorem}
    
  \newtheorem{remark}{\bf Remark}
 \newtheorem{lemma}{\bf Lemma}
  \newtheorem{prop}{\bf Proposition}
  \usepackage{dirtytalk}
\usepackage{amsmath}
\usepackage{amsfonts}
\usepackage{array}
\usepackage{mathtools,cuted}
\usepackage{amssymb}
\usepackage{bbm}
\usepackage{framed}
\usepackage{breqn}
\usepackage{cite}
\usepackage{enumitem}
\usepackage{diagbox}

\usepackage{cases}
\usepackage{tabu}
\usepackage{url}
\usepackage{multicol}
\usepackage{color}
\usepackage{pgfplots}
\usepackage[bottom]{footmisc}
\usepackage{environ}
\usepackage{subcaption}
\usepackage{nicefrac}
\usepackage{tikz}
\usetikzlibrary{arrows,decorations.markings}
\usetikzlibrary{decorations.pathreplacing}
\usetikzlibrary{automata,arrows,positioning,calc}
\title{On the Complexity of Sequential Incentive Design}
\author{$\text{Yagiz Savas}^{\star}$, Vijay Gupta, and Ufuk Topcu \thanks{ This work is supported in part by the grants DARPA D19AP00004 and ARL W911NF-17-2-0181. }\thanks{ We thank Melkior Ornik and Lillian Ratliff for their contributions to an earlier version of this paper. We thank Mustafa Karabag for giving the idea of the motion planning example and for insightful discussions.} \thanks{ Y. Savas and U. Topcu are with the Department of Aerospace Engineering, University of Texas at Austin, TX, USA. E-mail: \{yagiz.savas, utopcu\}@utexas.edu   } \thanks{V. Gupta is with the Department of Electrical Engineering, University of Notre Dame, IN, USA. E-mail: vgupta@nd.edu} }

\begin{document}
\maketitle
\begin{abstract}
    In many scenarios, a principal dynamically interacts with an agent and offers a sequence of incentives to align the agent's behavior with a desired objective. This paper focuses on the problem of synthesizing an incentive sequence that, once offered, induces the desired agent behavior even when the agent's intrinsic motivation is unknown to the principal. We model the agent's behavior as a Markov decision process, express its intrinsic motivation as a reward function, which belongs to a finite set of possible reward functions, and consider the incentives as additional rewards offered to the agent. We first show that the behavior modification problem (BMP), i.e., the problem of synthesizing an incentive sequence that induces a desired agent behavior at minimum total cost to the principal, is PSPACE-hard. Moreover, we show that by imposing certain restrictions on the incentive sequences available to the principal, one can obtain two NP-complete variants of the BMP. We also provide a sufficient condition on the set of possible reward functions under which the BMP can be solved via linear programming. Finally, we propose two algorithms to compute globally and locally optimal solutions to the NP-complete variants of the BMP. 
    
\end{abstract}
\section{Introduction}

Consider a scenario in which a principal offers \textit{a sequence of incentives} to influence the behavior of an agent with \textit{unknown} intrinsic motivation. For example, when an online retailer (principal) interacts with a customer (agent), the retailer may aim to convince the customer to purchase a number of products \textit{over time} by offering a sequence of discounts (incentives). The customer's willingness to purchase new products may depend on the ones purchased in the past, e.g., willingness to purchase a video game depends on whether the customer already owns a game console. Moreover, the retailer typically does not know what discount rates will encourage the customer to shop more. In this paper, we study the problem of synthesizing an incentive sequence that, once offered, induces the desired agent behavior \textit{regardless of the agent's intrinsic motivation while   minimizing the worst-case total cost to the principal}.

The design of incentives that align an agent's behavior with a principal's objective is a classical problem that has been studied from the perspective of control theory \cite{ho1982control}, game theory \cite{bacsar1984affine}, contract theory \cite{bolton2005contract}, and mechanism design \cite{nisan2001algorithmic} under various assumptions on related problem parameters. Despite the long history, to the best of our knowledge, there are only a few studies, e.g., \cite{conitzer2002complexity,conitzer2006computing,zhang2008value}, that analyze the computational complexity of incentive design problems. Such analyses are crucial to understand the computational challenges that naturally arise in incentive design problems and to develop effective algorithms for the synthesis of incentives. Here, we present a comprehensive complexity analysis of a particular class of sequential incentive design problems and propose several algorithms to solve them.

 We model the agent's behavior as a Markov decision process (MDP) \cite{Puterman}. MDPs model sequential decision-making under uncertainty and have been used to build recommendation systems \cite{shani2005mdp}, to develop sales promotion strategies \cite{rao1973dynamic}, and to design autonomous driving algorithms \cite{wei2011point}. By modelling the agent's behavior as an MDP, one can represent uncertain outcomes of the agent's actions. For instance, when a retailer makes a discount offer for purchases exceeding a certain amount, the actual amount the agent will spend upon accepting the offer can be expressed through a probability distribution.

 We express the agent's intrinsic motivation as a reward function and consider the incentives as additional nonnegative rewards offered for the agent's actions. Similar to the classical adverse selection problem \cite{bolton2005contract}, we assume that the agent's reward function is a private information and the principal only knows that the agent's reward function belongs to a \textit{finite} set of possible reward functions. The finite set assumption is standard in incentive design literature \cite{bergemann2005robust}; using certain techniques, e.g., clustering \cite{cadez2003model}, the principal may infer from historical data that agents can be categorized into a \textit{finite number of types} each of which is associated with a certain reward function. 


We consider a principal whose objective is to lead the agent to a desired set of target states in the MDP with maximum probability while minimizing its \textit{worst-case} expected total cost to induce such a behavior. In the considered setting, the worst-case scenario corresponds to a principal that interacts with an agent type that maximizes the incurred total cost over all possible agent types. The target set may represent, for example, the collection of products purchased by the customer. Then, the principal's objective corresponds to maximizing its profit while guaranteeing that the agent purchases the desired group of products regardless of its type.

This paper has four main contributions. First, we show that the behavior modification problem (BMP), i.e., the problem of synthesizing an incentive sequence that leads the agent to a desired target set with maximum probability while minimizing the worst-case total cost to the principal, is PSPACE-hard. Second, we show that by preventing the principal from adapting its incentive offers in time according to its history of interactions with the agent, one can obtain two ``easier" variants of the BMP which are NP-complete. Third, we show that, when the set of possible agent types contains a type that always demands an incentive offer that is higher than the ones demanded by any other type, the resulting BMP can be solved in polynomial-time via linear programming. Finally, we present two algorithms to solve the NP-complete variants of the BMP. The first algorithm computes a globally optimal solution to the considered problems based on a mixed-integer linear program. The second algorithm, on the other hand, computes a locally optimal solution by resorting to a variation \cite{lipp2016variations} of the so-called convex-concave procedure \cite{yuille2003concave} to solve a nonlinear optimization problem with bilinear constraints.

\noindent\textbf{Related Work:} Some of the results presented in Section \ref{poly_time_sec} of this paper has previously appeared in \cite{incentive_CDC} where we study a sequential incentive design problem in which the agent has a \textit{known} intrinsic motivation. In this paper, we present an extensive study of the sequential incentive design problem for agents with \textit{unknown} intrinsic motivation, which is a significantly different problem. In particular, we show that, when the agent's intrinsic motivation is unknown to the principal, it is, in general, PSPACE-hard to synthesize an optimal incentive sequence that induces the desired behavior. On the other hand, for an agent with known intrinsic motivation, an optimal incentive sequence can be synthesized in polynomial-time.

MDPs have been recently used to study the synthesis of sequential incentive offers in \cite{zhang2008value,reddi2012incentive,zhang2009policy,parkes2004mdp}. In \cite{zhang2008value}, the authors consider a principal with a \textit{limited budget} who aims to induce an agent behavior that maximizes the principal's utility. Similarly, in \cite{zhang2009policy}, a principal with a limited budget who aims to induce a \textit{specific} agent policy through incentive offers is considered. In \cite{reddi2012incentive}, the authors consider a multi-armed bandit setting in which a principal sequentially offers incentives to an agent with the objective of modifying its behavior. The paper \cite{parkes2004mdp} studies an online mechanism design problem in which a principal interacts with multiple agents and aims to maximize the expected social utility. Unlike the above mentioned studies, in this paper, we consider a principal \textit{with no budget constraints} that aims to modify the behavior of a single agent at minimum \textit{worst-case} total cost.

There are only a few results in the literature on the complexity of incentive design problems. In \cite{zhang2008value}, the authors consider a principal with budget constraints and prove, by a reduction from the Knapsack problem, that the considered problem is NP-hard even when the agent's reward function is known to the principal. Since we consider a principal with no budget constraints, the reduction techniques employed here are significantly different from the one used in \cite{zhang2008value}. In \cite{conitzer2006computing}, a number of NP-hardness results are presented for \textit{static} Stackelberg game settings in which the principal interacts with the agent \textit{only once}. In \cite{conitzer2002complexity}, the authors prove the complexity of several mechanism design problems in which the principal is allowed to condition its incentive offers on agent types. Even though the reduction techniques used in the above mentioned references are quite insightful, they cannot be applied to prove the complexity of the BMP which concerns a sequential setting in which the incentives are not conditioned on the agent types.

The complexity results presented in this paper are also closely related to the complexity of synthesizing optimal policies in robust MDPs with unstructured uncertainty sets \cite{le2007robust}, multi-model MDPs \cite{steimle2018multi}, and partially observable MDPs \cite{papadimitriou1987complexity}. In particular, the policy synthesis problems for all these models are PSPACE-hard. The main difference of the reduction technique we use in this paper from the ones used in the above mentioned references is the construction of the reward functions for each agent type. 
\section{Preliminaries}\label{prelim_sec}
 \noindent\textbf{Notation:} For a set $S$, we denote its cardinality by $\lvert S \rvert$. Additionally, $\mathbb{N}$$:=$$\{1,2,\ldots\}$, $\mathbb{R}$$:=$$(-\infty,\infty)$, and $\mathbb{R}_{\geq 0}$$:=$$[0,\infty)$. 
 \subsection{Markov decision processes}
{\setlength{\parindent}{0cm}
\begin{definition}
A \textit{Markov decision process} (MDP) is a tuple $\mathcal{M}$$:=$$(S, s_1, \mathcal{A}, \mathcal{P})$ where $S$ is a finite set of states, $s_1$$\in$$S$ is an initial state, $\mathcal{A}$ is a finite set of actions, and $ \mathcal{P}$$:$$S$$\times$$ \mathcal{A}$$\times$$S$$\rightarrow$$[0,1]$ is a transition function such that $\sum_{s'\in S}\mathcal{P}(s,a,s')$$=$$1$ for all $s$$\in$$S$ and $a$$\in$$\mathcal{A}$.
\end{definition}}
We denote the transition probability $ \mathcal{P}(s,a,s')$ by $ \mathcal{P}_{s,a,s'}$, and the set of available actions in $s$$\in$$S$ by $\mathcal{A}(s)$. A state $s$$\in$$S$ is \textit{absorbing} if $\mathcal{P}_{s,a,s}$$=$$1$ for all $a$$\in$$\mathcal{A}(s)$. The \textit{size of an MDP} is the number of triplets $(s,a,s')$$\in$$S$$\times$$\mathcal{A}$$\times$$S$ such that $\mathcal{P}_{s,a,s'}$$>$$0$.{\setlength{\parindent}{0cm}
\noindent \begin{definition}
For an MDP $\mathcal{M}$, a \textit{policy} is a sequence $\pi$$:=$$(d_1,d_2,d_3,\ldots)$ where each $d_t$$:$$ S$$\rightarrow$$\mathcal{A}$ is a \textit{decision rule} such that $d_t(s)$$\in$$\mathcal{A}(s)$ for all $s$$\in$$S$. A \textit{stationary} policy is of the form $\pi$$=$$(d,d,d,\ldots)$. We denote the set of all policies and all stationary policies by $\Pi(\mathcal{M})$ and $\Pi^S(\mathcal{M})$, respectively. 
\end{definition}}
We denote the action $a$$\in$$\mathcal{A}(s)$ taken by the agent in a state $s$$\in$$S$ under a stationary policy $\pi$ by $\pi(s)$.



\subsection{Incentive sequences} \label{incentive_prelim}
For an MDP $\mathcal{M}$, an \textit{information sequence} $I_t$ describes the information available to the principal at stage $t$$\in$$\mathbb{N}$ and is recursively defined as follows. At the first stage, an information sequence comprises only of the agent's current state, e.g., $I_1$$=$$(s_1)$. The principal makes incentive offers $\delta_1(I_1,a)$ to the agent for some actions, the agent takes an action $a_1$$\in$$\mathcal{A}(s_1)$ and transitions to a state $s_2$$\in$$S$. At the second stage, the information sequence becomes $I_2$$=$$(s_1,\gamma_1,a_1,s_2)$. The principal makes incentive offers $\delta_2(I_2,a)$ to the agent for some actions, the agent takes an action $a_2$$\in$$\mathcal{A}(s_2)$ and transitions to a state $s_3$$\in$$S$. The information sequence at stage $t$$\in$$\mathbb{N}$ is then recursively defined as $I_t$$=$$(s_1,\gamma_1,a_1,s_2,\ldots,s_{t-1},\gamma_{t-1},a_{t-1},s_t)$. We denote by $\mathcal{I}_t$ the set of all possible information sequences available to the principal at stage $t$$\in$$\mathbb{N}$.
{\setlength{\parindent}{0cm}
\noindent \begin{definition}
For an MDP $\mathcal{M}$, an \textit{incentive sequence} is a sequence $\gamma$$:=$$(\delta_1,\delta_2,\delta_3,\ldots)$ of incentive offers $\delta_t$$:$$\mathcal{I}_t$$\times $$\mathcal{A}$$\rightarrow$$\mathbb{R}_{\geq 0}$.  A \textit{stationary} incentive sequence is of the form $\gamma$$=$$(\delta,\delta,\delta,\ldots)$ where $\delta$$:$$S$$\times $$\mathcal{A}$$\rightarrow$$\mathbb{R}_{\geq 0}$. A \textit{stationary deterministic} incentive sequence is a stationary incentive sequence such that for all $s$$\in$$S$, $\sum_{a'\in\mathcal{A}(s)}\delta(s,a')$$=$$\delta(s,a)$ for some $a$$\in$$\mathcal{A}(s)$. For an MDP $\mathcal{M}$, we denote the set of all incentive sequences, all stationary incentive sequences, and all stationary deterministic incentive sequences by $\Gamma(\mathcal{M})$, $\Gamma^S(\mathcal{M})$, and $\Gamma^{SD}(\mathcal{M})$, respectively.
\end{definition}}
For a stationary incentive sequence $\gamma$$\in$$\Gamma^S(\mathcal{M})$, we denote the incentive offer for a state-action pair $(s,a)$ by $\gamma(s,a)$.

\subsection{Reachability}
An infinite sequence $\varrho^{\pi}$$=$$s_1s_2s_3\ldots$ of states generated in $\mathcal{M}$ under a policy $\pi$, which starts from the initial state $s_1$ and satisfies $\mathcal{P}_{s_t,d_t(s_t),s_{t+1}}$$>$$0$ for all $t$$\in$$\mathbb{N}$, is called a \textit{path}. We define the set of all paths in $\mathcal{M}$ generated under the policy $\pi$ by $Paths^{\pi}_{\mathcal{M}}$. We use the standard probability measure over the outcome set $Paths^{\pi}_{\mathcal{M}}$ \cite{Model_checking}. Let $\varrho^{\pi}[t]$$:=$$s_t$ denote the state visited at the $t$-th stage along the path $\varrho^{\pi}$. We define 
\begin{align*}
     \text{Pr}_{\mathcal{M}}^{\pi}(Reach[B]):=\text{Pr}\{\varrho^{\pi}\in Paths^{\pi}_{\mathcal{M}}: \exists t\in \mathbb{N}, \varrho^{\pi}[t]\in B \}
\end{align*} 
as the probability with which the paths generated in $\mathcal{M}$ under $\pi$ reaches the set $B$$\subseteq$$S$. We denote the maximum probability of reaching the set $B$ under any policy $\pi$$\in$$\Pi(\mathcal{M})$ by
\begin{align*}
    R_{\max}(\mathcal{M},B):=\max_{\pi\in\Pi(\mathcal{M})} \text{Pr}_{\mathcal{M}}^{\pi}(Reach[B]).
\end{align*}

The existence of maximum in the above definition follows from Lemma 10.102 in \cite{Model_checking}, and the value of $R_{\max}(\mathcal{M},B)$ can be computed via linear programming \cite[Chapter 10]{Model_checking}. 

\section{Problem Statement}\label{sec_prob_formulate}

We consider an \textit{agent} whose behavior is modeled as an MDP $\mathcal{M}$ and a \textit{principal} that provides a sequence $\gamma$$\in$$\Gamma(\mathcal{M})$ of incentives to the agent in order to induce a desired behavior.


We assume that the agent's intrinsic motivation is \textit{unknown} to the principal, but it can be expressed by one of finitely many reward functions. Formally, let $\Theta$ be a finite set of agent \textit{types}, and for each $\theta$$\in$$\Theta$, $\mathcal{R}_{\theta}$$:$$S$$\times$$\mathcal{A}$$\rightarrow$$\mathbb{R}$ be the associated reward function. The principal knows the function $\mathcal{R}_{\theta}$ associated with each type $\theta$$\in$$\Theta$, but does not know the \textit{true agent type} $\theta^{\star}$$\in$$\Theta$.

The sequence of interactions between the principal and the agent is as follows. At stage $t$$\in$$\mathbb{N}$, the agent occupies a state $s_t$$\in$$S$. The information sequence $I_t$, defined in Section \ref{incentive_prelim}, describes the information available to the principal about the agent. The principal offers an incentive $\delta_t(I_t,a)$ to the agent for each action $a$$\in$$\mathcal{A}$, and the agent chooses an action  $a_t$$\in$$\mathcal{A}(s_t)$ to maximize its immediate total reward, i.e.,
\begin{align}\label{myopic_objective}
    a_t\in \arg\max_{a\in\mathcal{A}(s_t)} \Big[\mathcal{R}_{\theta^{\star}}(s_t,a)+\delta_t(I_t,a)\Big].
\end{align}
Depending on the action taken by the agent, the principal pays $\delta_t(I_t,a_t)$ to the agent. Finally, the agent transitions to the next state $s_{t+1}$$\in$$S$ with probability $\mathcal{P}_{s_t,a_t,s_{t+1}}$, the principal updates its information sequence to $I_{t+1}$$=$$(I_t,\delta_t,a_t,s_{t+1})$, and so on. 

The principal offers a sequence of incentives to the agent so that the agent reaches a target state with maximum probability \textit{regardless of its type} and the worst-case expected total cost of inducing the desired behavior to the principal is minimized. 

{\setlength{\parindent}{0cm}
\begin{problem} \textbf{(Behavior modification problem (BMP))} For an MDP $\mathcal{M}$, a set $B$$\subseteq$$S$ of absorbing target states, and a set $\Theta$ of possible agent types, synthesize an incentive sequence $\gamma$$\in$$\Gamma(\mathcal{M})$ that leads the agent to a target state with maximum probability while minimizing the expected total cost, i.e.,
\begin{subequations}
\begin{flalign}\label{reach_obj}
    \underset{\gamma\in\Gamma(\mathcal{M})}{\text{minimize}}\max_{\theta\in\Theta}&\   \mathbb{E}^{\pi^{\star}}\Bigg[\sum_{t=1}^{\infty}\delta_t(I_t,A_t)\Big| \theta\Bigg]&&\\
    \text{subject to:} &\ \label{policy_opt_def} \pi^{\star}=(d_1^{\star},d_2^{\star},d_3^{\star},\ldots)&&\\ \label{opt_set_agent}
   \forall t\in\mathbb{N},\forall s\in S, &   \   d_t^{\star}(s)\in \arg\max_{a\in\mathcal{A}(s)} \Big[\mathcal{R}_{\theta}(s,a)+\delta_t(I_t,a)\Big]&&\raisetag{18pt}\\ \label{last_cons}
    & \  \text{Pr}^{\pi^{\star}}_{\mathcal{M}}(Reach[B])=R_{\max}(\mathcal{M},B).&&
\end{flalign}
\end{subequations}
\end{problem}}
The expectation in \eqref{reach_obj} is taken over the paths that are induced by the policy $\pi^{\star}$, i.e., the principal pays the offered incentive if and only if the agent takes the incentivized action. 

{\setlength{\parindent}{0cm}
\begin{remark} The decision-making model \eqref{myopic_objective} describes a \textit{myopic} agent which aims to maximize only its immediate rewards. Clearly, this model is restrictive in terms of its capability of explaining sophisticated agent behaviors. The main reason that we use such a simple model is to avoid notational burden that comes with more expressive alternatives. In Appendix \ref{Appendix_remarks}, we show that the behavior modification of an agent with a finite decision horizon is \textit{computationally not easier} than the behavior modification of a myopic agent. We also show how the solution techniques developed for the behavior modification of a myopic agent can be utilized to modify the behavior of an agent with a finite decision horizon.
\end{remark}}

{\setlength{\parindent}{0cm}
\begin{remark} The BMP requires all target states $s$$\in$$B$ to be absorbing. We impose such a restriction just to avoid the additional definitions and notations required to consider non-absorbing target states. In Appendix \ref{Appendix_remarks}, we show how the results presented in this paper can be applied to MDPs with non-absorbing target states.
\end{remark}}

{\setlength{\parindent}{0cm}
\begin{remark} The BMP inherently assumes that the principal knows the transition function $\mathcal{P}$ of the given MDP $\mathcal{M}$. It is clear that, by allowing $\mathcal{P}$ to belong to an uncertainty set, one cannot obtain a computationally easier problem. 
\end{remark}}

\section{Summary of the Results}

\begin{table*}\centering
\caption{A summary of the presented results.  }
\begin{tabular}{ |>{\centering\arraybackslash}p{7.2 cm} || c|| c|| c| } 
 \hline
 {{\textbf{Problem}}} & \textbf{Complexity} & \textbf{Globally optimal solution}&  \textbf{Locally optimal solution} \\ 
 \hline
   \hline 
 Behavior modification (BMP) & PSPACE-hard & --- & ---\\ 
 \hline
 Non-adaptive behavior modification (N-BMP)& NP-complete & MILP & Convex-concave procedure \\ 
  \hline
 Non-adaptive single-action behavior modification (NS-BMP)& NP-complete & MILP & Convex-concave procedure \\
  \hline
 Behavior modification of a dominant type (BMP-D) & P & LP & --- \\
 \hline
\end{tabular}
\label{table_results}
\end{table*}

In this section, we provide a brief summary of the presented results. We keep the exposition non-technical to improve readability; precise statements are provided in later sections. Table \ref{table_results} illustrates an overview of the considered problems, their complexity, and the proposed solution techniques.

We first show that the BMP is PSPACE-hard. That is, unless P=PSPACE, an event considered to be less likely than the event P=NP \cite{papadimitriou1987complexity}, an efficient algorithm to solve the BMP does not exist. Due to this discouraging result, instead of studying the BMP in its most generality, we focus on its variants in which we restrict the incentive sequences available to the principal.

As the first restriction, we allow the principal to use only stationary incentive sequences to modify the agent's behavior and refer to the resulting problem as the \textit{non-adaptive behavior modification problem (N-BMP)}. We show that the decision problem associated with the N-BMP is NP-complete \textit{even when the given MDP has only deterministic transitions}. To find a globally optimal solution to the N-BMP, we formulate a mixed-integer linear program (MILP) in which the integer variables correspond to actions taken by different agent types in a given state. We also show that the N-BMP can be formulated as a nonlinear optimization problem with bilinear constraints for which a locally optimal solution can be obtained using the so-called convex-concave procedure \cite{lipp2016variations}. 

As the second restriction, we allow the principal to use only stationary deterministic incentive sequences to modify the agent's behavior and refer to the resulting problem as the \textit{non-adaptive single-action behavior modification problem (NS-BMP)}. 
We prove that the decision problem associated with the NS-BMP is NP-complete \textit{even when the agent has state-independent reward functions}. We also show that globally and locally optimal solutions to the NS-BMP can be computed by slightly modifying the methods developed to solve the N-BMP.

Finally, we consider the case in which the set of agent types include a dominant type which always demands the principal to offer an incentive amount that is higher than the ones demanded by any other agent type. We prove that solving the BMP instances that involve a dominant type is equivalent to modifying the behavior of the dominant type. We show that the behavior modification problem of a dominant type (BMP-D) is in P, and present an approach based on a linear program (LP) to solve the BMP-D.
\section{Optimal Agent Behavior and Optimal Incentive Sequences}
In the BMP, we aim to synthesize an incentive sequence that induces an optimal agent policy that reaches a target set with maximum probability while minimizing the expected total cost to the principal. To obtain a well-defined BMP, in this section, we first precisely specify how the agent behaves when there are multiple optimal policies. We then show that, in general, an incentive sequence that minimizes the expected total cost to the principal may not exist. Hence, we focus on $\epsilon$-optimal incentive sequences where $\epsilon$$>$$0$ is an arbitrarily small constant.

\subsection{Agent's behavior when multiple optimal policies exist}
\label{multi_policy_sec} For a given incentive sequence, the agent's optimal policy may not be unique. In the presence of multiple optimal policies, there are only two possible cases: either (i) there exists an optimal policy that violates the \textit{reachability constraint} in \eqref{last_cons} or (ii) all optimal policies satisfy the constraint in \eqref{last_cons}.

\begin{figure}[t!]
\begin{subfigure}[t]{0.33\linewidth}\vspace{-2.15cm}
\scalebox{0.9}{
\begin{tikzpicture}[->, >=stealth', auto, semithick, node distance=2cm]

    \tikzstyle{every state}=[fill=white,draw=black,thick,text=black,scale=0.7]

    \node[state,initial,initial text=] (s_1) {$s_1$};
    \node[state] (s_2) [right=17mm of s_1]  {$s_2$};

\path
(s_1)  edge  [loop above=10]    node{$a_1, 0$}     (s_1)
(s_1)	 edge     node{$a_2, -1$}     (s_2)
(s_2)  edge  [loop right=10]    node{}     (s_2);
\end{tikzpicture}}
\end{subfigure}
\hspace{0.13\linewidth}
\begin{subfigure}[t]{0.33\linewidth}
\scalebox{0.9}{
\begin{tikzpicture}[->, >=stealth', auto, semithick, node distance=2cm]

    \tikzstyle{every state}=[fill=white,draw=black,thick,text=black,scale=0.7]

    \node[state,initial,initial text=] (s_1) {$s_1$};
    \node[state] (s_2) [right=17mm of s_1]  {$s_2$};

\path
(s_1)  edge  [loop above=10]    node{$a_1, 0$}     (s_1)
(s_1)	 edge     node{$a_2, -1$}     (s_2)
(s_1)  edge[dashed, loop below=10]    node{$a_3, -1$}     (s_1)
(s_1)	 edge[dashed, out=290, in=200]     node{}     (s_2)
(s_2)  edge  [out=270,loop right=10]    node{}     (s_2);
\end{tikzpicture}}
\end{subfigure}
\caption{MDP examples to illustrate the agent's behavior in the existence of multiple optimal policies. Solid lines represent deterministic transitions; dashed lines represent transitions with equal probability. The initial state is $s_1$. The type set $\Theta$ satisfies $\lvert \Theta \rvert$$=$$1$, i.e., the principal knows the true agent type. The target set is $B$$=$$\{s_2\}$. The tuples $(a,r)$ next to the arrows indicate the action $a$ and the reward $r$.}
\label{fig:multiple_policies}
\end{figure}
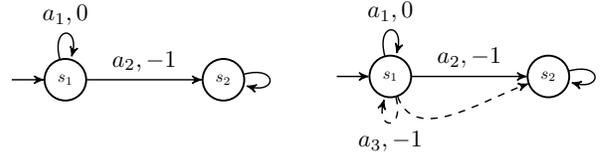

We first analyze case (i). Consider the MDP given in Fig. \ref{fig:multiple_policies} (left). Let the agent's reward function be $\mathcal{R}_{\theta^{\star}}(s_1,a_1)$$=$$0$ and $\mathcal{R}_{\theta^{\star}}(s_1,a_2)$$=$$-1$. That is, in the absence of incentives, it is optimal for the agent to stay in state $s_1$ under action $a_1$. 

Suppose that the principal offers the \textit{stationary} incentives $\gamma(s_1,a_1)$$=$$0$ and $\gamma(s_1,a_2)$$=$$1$ to the agent. Then, we have
\begin{align*}
    \{a_1,a_2\}=\arg\max_{a\in\mathcal{A}(s_1)} \Big[\mathcal{R}_{\theta^{\star}}(s_1,a)+\gamma(s_1,a)\Big]
\end{align*}
which implies that the agent has multiple optimal policies. Note that under the optimal stationary policy $\pi(s_1)$$=$$a_2$, the agent reaches the target state $s_2$ with probability 1, whereas under the optimal stationary policy $\pi(s_1)$$=$$a_1$, it reaches the target state $s_2$ with probability 0.
We assume that, in such a scenario, the agent behaves adversarially against the principal.

\noindent \textbf{Assumption:} Under the provided incentive sequence, if there exists an optimal policy following which the agent can violate the constraint in \eqref{last_cons}, then the agent follows such a policy. 

We now analyze case (ii). Consider the MDP given in Fig. \ref{fig:multiple_policies} (right). In this MDP, in addition to the actions $a_1$ and $a_2$, the agent can take a third action $a_3$ which leads the agent to the states $s_1$ and $s_2$ with equal probability. Let the reward function be $\mathcal{R}_{\theta^{\star}}(s_1,a_1)$$=$$0$, $\mathcal{R}_{\theta^{\star}}(s_1,a_2)$$=$$-1$, and $\mathcal{R}_{\theta^{\star}}(s_1,a_3)$$=$$-1$.

Suppose that the principal offers the stationary incentives $\gamma(s_1,a_1)$$=$$0$, $\gamma(s_1,a_2)$$=$$2$, and $\gamma(s_1,a_3)$$=$$2$. Then, we have
\begin{align*}
    \{a_2,a_3\}=\arg\max_{a\in\mathcal{A}(s_1)} \Big[\mathcal{R}_{\theta^{\star}}(s_1,a)+\gamma(s_1,a)\Big],
\end{align*}
and the agent has multiple optimal policies. Under all its optimal policies, the agent reaches the target state $s_2$ with probability 1. However, if the agent follows the stationary policy $\pi(s_1)$$=$$a_2$, the expected total cost to the principal is equal to 2, whereas the stationary policy $\pi(s_1)$$=$$a_3$ incurs the expected total cost of 4. In such a scenario, we again assume that the agent behaves adversarially against the principal. 

\noindent \textbf{Assumption:} Under the provided incentive sequence, if all optimal agent policies satisfy the constraint in \eqref{last_cons}, then the agent follows the policy that maximizes the expected total cost to the principal.

\subsection{Non-existence of optimal incentive sequences}\label{non_existence_opt_section}
We now illustrate with an example that, in general, there may exist no incentive sequence $\gamma$$\in$$\Gamma(\mathcal{M})$ that attains the minimum in \eqref{reach_obj}-\eqref{last_cons}. Consider again the MDP given in Fig. \ref{fig:multiple_policies} (left). In this example, without loss of generality, we can focus on stationary incentive sequences to minimize the expected total cost to the principal. The feasible set of stationary incentive sequences is described by the set 
\begin{align*}
    \{\gamma(s,a)\geq 0 \ |\  \gamma(s_1,a_2)-\gamma(s_1,a_1)\geq 1+\tilde{\epsilon},\  \tilde{\epsilon}>0\}
\end{align*}
which is not a compact set due to the condition $\tilde{\epsilon}$$>$$0$. The incentive offers satisfying $\gamma(s_1,a_2)$$-$$\gamma(s_1,a_1)$$=$$1$, i.e., $\tilde{\epsilon}$$=$$0$, do not belong to the feasible set because under such incentive sequences the agent adversarially follows the stationary policy $\pi(s_1)$$=$$a_1$ which violates the constraint in \eqref{last_cons}.

Motivated by the above example, in the rest of the paper, we focus on obtaining an $\epsilon$-optimal solution to the BMP. 
{\setlength{\parindent}{0cm}
\begin{definition}
Let $\mathcal{D}$ be a non-empty set, $h$$:$$\mathcal{D}$$\rightarrow$$\mathbb{R}_{\geq 0}$ be a function, and $\epsilon$$>$$0$ be a constant. A point $\overline{x}$$\in$$\mathcal{D}$ is said to be an \textit{$\epsilon$-optimal solution} to the problem $\min_{x\in\mathcal{D}}h(x)$ if 
\begin{align*}
    h(\overline{x})\leq \min_{x\in\mathcal{D}}h(x)+\epsilon.
\end{align*}
\end{definition}}
\section{Complexity of Behavior Modification}
In this section, we present the results on the computational complexity of the BMP and its two variants. Before proceeding with the complexity results, we first show that a feasible solution to the BMP can be computed efficiently.

\subsection{A feasible solution to the behavior modification problem}\label{feasible_sectiob}
Let $C$$:$$S$$\times$$\mathcal{A}$$\rightarrow$$\mathbb{R}_{\geq 0}$ be defined as
\begin{align}\label{cost_cost_cost}
    C(s,a):=\max_{\theta\in\Theta}\Big(\max_{a'\in\mathcal{A}(s)}\mathcal{R}_{\theta}(s,a')-\mathcal{R}_{\theta}(s,a)+\epsilon\Big).
\end{align}
The value of $C(s,a)$ denotes an incentive amount under which the agent takes the incentivized action $a$$\in$$\mathcal{A}(s)$ in state $s$$\in$$S$ regardless of its type. To see this, note that the term inside the parentheses denotes an amount that, once offered, makes the action $a$$\in$$\mathcal{A}(s)$ uniquely optimal for the agent type $\theta$$\in$$\Theta$. Hence, by taking the maximum over all agent types, we ensure that, under the stationary incentive offer $\gamma(s,a)$$=$$C(s,a)$, the agent takes the action $a$$\in$$\mathcal{A}(s)$ regardless of its type.

We can obtain a feasible solution to the BMP as follows. First, compute a policy $\pi^{\star}$$\in$$\Pi(\mathcal{M})$ that reaches the desired target set $B$ with maximum probability $R_{\max}(\mathcal{M},B)$ by solving a linear program \cite[Chapter 10]{Model_checking}. Then, offer the stationary incentive sequence 
\begin{align}
    \gamma(s,a)=\begin{cases} C(s,a) & \text{if}\  \pi^{\star}(s)=a\\
    0 & \text{otherwise}
    \end{cases}
\end{align}
to the agent. Under the above incentive sequence, the agent follows the desired policy $\pi^{\star}$ regardless of its type; hence, it reaches the target set $B$ with maximum probability. 

There are two important aspects of the incentive sequence synthesized above. First, it is computed in time polynomial in the size of $\mathcal{M}$ and $\Theta$. Second, it conservatively incentivizes the same actions for all agent types ignoring the principal's knowledge on the reward functions $\mathcal{R}_{\theta}$. In the next section, we show that, in order to exploit the principal's knowledge on the reward functions and minimize the incurred total cost, one should solve a more challenging computational problem.

\subsection{Complexity results}

We first present a discouraging result which shows that, unless P=PSPACE, solving the BMP in time polynomial in the size of $\mathcal{M}$ and $\Theta$ is not possible.
{\setlength{\parindent}{0cm} 
\begin{thm} \label{PSPACE_thm} For $N$$\in$$\mathbb{N}$, deciding whether there exists a feasible solution to the BMP with the objective value less than or equal to $N$ is PSPACE-hard.
\end{thm}}

We provide a proof for the above result in Appendix \ref{appendix_proofs} where we show that the PSPACE-complete quantified satisfiability problem (QSAT) \cite{sipser1996introduction} can be reduced to the BAP. The reduction closely follows the reduction of QSAT to \textit{partially observable MDP} problems \cite{papadimitriou1987complexity} in which the objective is to synthesize an observation-based controller to maximize the expected total reward collected by an agent. 


Next, we introduce two variants of the BMP in which the principal is allowed to use only a subset of the set $\Gamma(\mathcal{M})$ of all incentive sequences.
{\setlength{\parindent}{0cm}
\begin{definition}
When the principal is restricted to use stationary incentive sequences, i.e., the minimization in \eqref{reach_obj} is performed over the set $\Gamma^S(\mathcal{M})$, the resulting BMP is referred to as the \textit{non-adaptive BMP} (N-BMP).
\end{definition}}

In N-BMP, the principal does not utilize its past experiences, i.e., information sequence $I_t$, with the agent to refine the incentive offers in the future. Therefore, in a sense, the N-BMP describes an offline incentive design problem. The following result shows that the N-BMP is computationally ``easier" than the BMP. However, it is still not likely to be solved efficiently even for MDPs with deterministic transition functions.
{\setlength{\parindent}{0cm}
\begin{thm} \label{N-BMP_thm} For $N$$\in$$\mathbb{N}$, deciding whether there exists a feasible solution to the N-BMP with the objective value less than or equal to $N$ is NP-complete even when the transition function $\mathcal{P}$ satisfies $\mathcal{P}_{s,a,s'}$$\in$$\{0,1\}$ for all $s,s'$$\in$$S$ and all $a$$\in$$\mathcal{A}$.
\end{thm}}
A proof for the above result is provided in Appendix \ref{appendix_proofs} where we reduce the NP-complete Euclidean path traveling salesman problem \cite{papadimitriou1977euclidean} to the N-BMP. The reduction demonstrates that, to lead the agent to a target state at minimum total cost, the principal may need to reveal the agent's true type.  


{\setlength{\parindent}{0cm}
\begin{definition}
When the principal is restricted to use stationary deterministic incentive sequences, i.e., minimization in \eqref{reach_obj} is performed over the set $\Gamma^{SD}(\mathcal{M})$, the resulting BMP is referred to as the \textit{non-adaptive single-action BMP} (NS-BMP).
\end{definition}}

In NS-BMP, the principal offers incentives \textit{only for a single action} in a given state. Such a restriction can be encountered in practice, for example, when a retailer offers discounts only for the purchases that exceed a certain total sum. The following result shows that the NS-BMP is a challenging problem even when the agent has a \textit{state-independent} reward function.
{\setlength{\parindent}{0cm}
\begin{thm}\label{NS-BMP_thms} For $N$$\in$$\mathbb{N}$, deciding whether there exists a feasible solution to the NS-BMP with the objective value less than or equal to $N$ is NP-complete even when $\mathcal{R}_{\theta}(s,a)$$=$$\mathcal{R}_{\theta}(t,a)$ for all $s,t$$\in$$S$, all $a$$\in$$\mathcal{A}$, and all $\theta$$\in$$\Theta$.
\end{thm}}
A proof for the above result is provided in Appendix \ref{appendix_proofs} where we reduce the NP-complete set cover problem \cite{papadimitriou1987complexity} to the NS-BMP. The main idea in the reduction is that, to lead the agent to a target state, the principal may need to reveal the true agent type \textit{by the end of a fixed number of stages}.

\section{A Sufficient Condition for Polynomial-Time Solvability} \label{poly_time_sec}
In this section, we derive a sufficient condition on the structure of the type set $\Theta$ under which the BMP can be solved by assuming that the true agent type is known to the principal. We define the BMP of a dominant type (BMP-D) as the BMP instances that satisfy the derived sufficient condition and show that an $\epsilon$-optimal solution to the BMP-D can be computed in time polynomial in the size of $\mathcal{M}$ and $\Theta$.

\subsection{Behavior modification of a dominant type: formulation }
Let $f$$:$$\Gamma(\mathcal{M})$$\times$$\Theta$$\rightarrow$$\mathbb{R}$ be a function such that
\begin{align}\label{cost_function}
    f(\gamma,\theta):=\mathbb{E}^{\pi^{\star}}\Bigg[\sum_{t=1}^{\infty}\delta_t(I_t,A_t)\Big| \theta\Bigg]
\end{align}
is the expected total cost incurred by the principal when the sequence of incentive offers is $\gamma$$\in$$\Gamma(\mathcal{M})$ and the agent type is $\theta$$\in$$\Theta$. In \eqref{cost_function}, $\pi^{\star}$$\in$$\Pi(\mathcal{M})$ is the agent's optimal policy under the incentive sequence $\gamma$. Moreover, let $\Gamma_{R,\theta}(\mathcal{M})$$\subseteq$$\Gamma(\mathcal{M})$ be the set of incentive sequences under which the agent type $\theta$$\in$$\Theta$ reaches the target set $B$$\subseteq$$S$ with maximum probability. Precisely, we have $\gamma$$\in$$\Gamma_{R,\theta}(\mathcal{M})$ if and only if the optimal policy $\pi^{\star}$ of the agent type $\theta$ under the sequence $\gamma$$\in$$\Gamma(\mathcal{M})$ satisfies the equality in \eqref{last_cons}. Finally, let 
\begin{align}
\Gamma_{R,\Theta}(\mathcal{M}):=\bigcap\limits_{\theta\in\Theta}\Gamma_{R,\theta}(\mathcal{M})
\end{align}
be the set of incentive sequences under which the agent reaches the target set with maximum probability \textit{regardless of its type}.
Then, the BMP in \eqref{reach_obj}-\eqref{last_cons} can equivalently be written as 
\begin{align}
    \min_{\gamma\in\Gamma_{R,\Theta}(\mathcal{M})}\max_{\theta\in\Theta}f(\gamma,\theta).
\end{align}

Using weak duality \cite{boyd2004convex}, together with the fact that $\Gamma_{R,\Theta}$$\subseteq$$\Gamma_{R,\theta}$ for all $\theta$$\in$$\Theta$, we have 
\begin{align}\label{weak_duality}
    \min_{\gamma\in\Gamma_{R,\Theta}(\mathcal{M})}\max_{\theta\in\Theta}f(\gamma,\theta)\geq \max_{\theta\in\Theta}\min_{\gamma\in\Gamma_{R,\theta}(\mathcal{M})}f(\gamma,\theta).
\end{align}
The above inequality implies that a principal that \textit{knows} the true agent type $\theta^{\star}$ can always modify the agent's behavior at an expected total cost that is less than or equal to the worst-case expected total cost incurred by a principal that \textit{does not know} the true agent type. We now provide a sufficient condition under which the knowledge of the true agent type $\theta^{\star}$ does not help the principal to strictly decrease its expected total cost, i.e., the inequality in \eqref{weak_duality} becomes an equality. For each $\theta$$\in$$\Theta$ and $s$$\in$$S$, let $\mathcal{R}^{\max}_{\theta}(s)$$:=$$\max_{a'\in\mathcal{A}(s)}\mathcal{R}_{\theta}(s,a')$.
{\setlength{\parindent}{0cm}
\begin{thm}\label{theorem_dominant}
For a given MDP $\mathcal{M}$ and a type set $\Theta$, if there exists $\theta_d$$\in$$\Theta$ such that, 
\begin{align}
    \mathcal{R}^{\max}_{\theta_d}(s)- \mathcal{R}_{\theta_d}(s,a) \geq     \mathcal{R}^{\max}_{\theta}(s)- \mathcal{R}_{\theta}(s,a)\label{condition_dominant}
\end{align}
for all $\theta$$\in$$\Theta$, $s$$\in$$S$, and $a$$\in$$\mathcal{A}$, then $\Gamma_{R,\theta_d}(\mathcal{M})$$=$$\Gamma_{R,\Theta}(\mathcal{M})$ and
\begin{align*}
    \min_{\gamma\in\Gamma_{R,\Theta}(\mathcal{M})}\max_{\theta\in\Theta}f(\gamma,\theta)= \min_{\gamma\in\Gamma_{R,\theta_d}(\mathcal{M})}f(\gamma,\theta_d).
\end{align*}
\end{thm}}
{\setlength{\parindent}{0cm}
\begin{definition}\label{dominant_problem_def} When there exists a \textit{dominant type} $\theta_d$$\in$$\Theta$ that satisfies the condition in \eqref{condition_dominant}, the resulting BMP is referred to as the \textit{BMP of a dominant type} (BMP-D).
\end{definition}}

Theorem \ref{theorem_dominant} implies that one can solve the BMP instances that satisfy the condition in \eqref{condition_dominant} by assuming that the true agent type $\theta^{\star}$ is the type $\theta_d$. Hence, solving the BMP-D is equivalent to finding an incentive sequence $\gamma^{\star}$$\in$$\Gamma(\mathcal{M})$ such that 
\begin{align}\label{dominant_problem_opt}
   \gamma^{\star}\in\arg\min_{\gamma\in\Gamma_{R,\theta_d}(\mathcal{M})}f(\gamma,\theta_d).
\end{align}
\subsection{Behavior modification of a dominant type: reformulation}\label{reformulation_section}
In general, the minimum in \eqref{dominant_problem_opt} is not attainable by any feasible incentive sequence. Hence, we aim to find an $\epsilon$-optimal solution as defined as in Section \ref{non_existence_opt_section}. In this section, we show that finding an $\epsilon$-optimal solution to the BMP-D is equivalent to solving a certain constrained MDP problem.

For a given MDP $\mathcal{M}$, we first partition the set $S$ of states into three disjoint sets as follows. Let $B$$\subseteq$$S$ be the set of target states and $S_0$$\subseteq$$S$ be the set of states that have zero probability of reaching the states in $B$ under any policy. Finally, we let $S_r$$=$$S\backslash B\cup S_0$. These sets can be found in time polynomial in the size of the MDP using graph search algorithms \cite{Model_checking}.

For $\overline{\epsilon}$$\geq$$0$, let $\phi_{\overline{\epsilon}}$$:$$S\times \mathcal{A}$$\rightarrow$$\mathbb{R}_{\geq0}$ be \textit{the cost of control function} 
\begin{align*}
   \phi_{\overline{\epsilon}}(s,a):=
   \begin{cases} \mathcal{R}_{\theta_d}^{\max}(s)-\mathcal{R}_{\theta_d}(s,a)+\overline{\epsilon}& \text{if}\ s\in S_r, \ a\in\mathcal{A}(s)\\
   0 & \text{otherwise}.
    \end{cases}
\end{align*}
Additionally, let $\Pi_{R,\theta_d}(\mathcal{M})$ be the set of policies under which the dominant agent type $\theta_d$ reaches the target set $B$ with maximum probability, i.e., $\pi'$$\in$$\Pi_{R,\theta_d}(\mathcal{M})$ if and only if 
\begin{align*}
    \text{Pr}^{\pi'}_{\mathcal{M}}(Reach[B])=R_{\max}(\mathcal{M},B).
\end{align*}

For any incentive sequence $\gamma$$\in$$\Gamma(\mathcal{M})$, there exists a corresponding optimal policy $\pi$$\in$$\Pi(\mathcal{M})$ for the agent type $\theta_d$. Moreover, for the principal to induce an agent policy $\pi$$=$$(d_1,d_2,\ldots)$ such that $\pi$$\in$$\Pi_{R,\theta_d}(\mathcal{M})$, it is \textit{necessary} that, for any information sequence $I_t$$=$$(s_1,\gamma_1,a_1,\ldots,s_t)$, we have $\delta_t(I_t,a)$$\geq$$\phi_0(s_t,a)$ if $d_t(s_t)$$=$$a$. Otherwise, we have
\begin{align*}
   \mathcal{R}_{\theta_d}(s_t,a)+\delta_t(I_t,a)< \mathcal{R}^{\max}_{\theta_d}(s_t),
\end{align*}
 which implies that the action $a$$\in$$\mathcal{A}(s_t)$ cannot be optimal for the agent in state $s_t$$\in$$S$ at stage $t$$\in$$\mathbb{N}$. Let $g_{\overline{\epsilon}}$$:$$\Pi(\mathcal{M})$$\rightarrow$$\mathbb{R}_{\geq 0}$ be \textit{the expected total cost of control function} 
\begin{align}\label{total_cost_function}
    g_{\overline{\epsilon}}(\pi):=\mathbb{E}^{\pi}\Bigg[\sum_{t=1}^{\infty}\phi_{\overline{\epsilon}}(S_t,A_t)\Big| \theta_d\Bigg]
\end{align}
where the expectation is taken over the paths induced by the dominant agent type's policy $\pi$.
The necessary condition described above implies that we have
\begin{align}\label{necessary_ineq}
    \min_{\gamma\in\Gamma_{R,\theta_d}(\mathcal{M})}f(\gamma,\theta_d)\geq \min_{\pi\in\Pi_{R,\theta_d}(\mathcal{M})}g_0(\pi).
\end{align}

Similarly, for the principal to induce an agent policy $\pi$$=$$(d_1,d_2,\ldots)$ such that $\pi$$\in$$\Pi_{R,\theta_d}(\mathcal{M})$, it is \textit{sufficient} that, for any information sequence $I_t$$=$$(s_1,\gamma_1,a_1,\ldots,s_t)$ where $s_t$$\in$$S_r$, we have $\delta_t(I_t,a)$$=$$\phi_{\overline{\epsilon}}(s_t,a)$ for some $\overline{\epsilon}$$>$$0$ if $d_t(s_t)$$=$$a$. This follows from the fact that, if $\delta_t(I_t,a)$$=$$\phi_{\overline{\epsilon}}(s_t,a)$, then 
\begin{align*}
   \mathcal{R}_{\theta_d}(s_t,a)+\delta_t(I_t,a)= \mathcal{R}^{\max}_{\theta_d}(s_t)+\overline{\epsilon}>\mathcal{R}_{\theta_d}(s_t,a')
\end{align*}
for all $a'$$\in$$\mathcal{A}(s)\backslash \{a\}$.
Hence, the action $a$$\in$$\mathcal{A}(s_t)$ is uniquely optimal for the agent in state $s_t$$\in$$S_r$ at stage $t$$\in$$\mathbb{N}$. The sufficient condition described above implies that, for any $\overline{\epsilon}$$>$$0$, 
\begin{align}\label{sufficient_ineq}
  \min_{\pi\in\Pi_{R,\theta_d}(\mathcal{M})}g_{\overline{\epsilon}}(\pi)\geq  \min_{\gamma\in\Gamma_{R,\theta_d}(\mathcal{M})}f(\gamma,\theta_d).
\end{align}

We now show that finding an $\epsilon$-optimal solution to the BMP-D is equivalent to solving the optimization problem on the left hand side of \eqref{sufficient_ineq} for some $\overline{\epsilon}$$>$$0$. Let $\Pi^S_{R,\theta_d}(\mathcal{M})$$\subseteq$$\Pi_{R,\theta_d}(\mathcal{M})$ be the set of \textit{stationary} policies under which the dominant agent type $\theta_d$ reaches the target set with maximum probability. 
 {\setlength{\parindent}{0cm}
\begin{lemma}\cite{stochastic_safest}\label{stationary sufficient}
For any $\overline{\epsilon}$$\geq$$0$, 
\begin{align}
    \min_{\pi\in\Pi_{R,\theta_d}(\mathcal{M})}g_{\overline{\epsilon}}(\pi)=\min_{\pi\in\Pi^S_{R,\theta_d}(\mathcal{M})}g_{\overline{\epsilon}}(\pi).
\end{align}
\end{lemma}}
The above result, which is proven in \cite{stochastic_safest}, implies that stationary policies are sufficient to minimize the expected total cost function \eqref{total_cost_function} subject to the reachability constraint in \eqref{last_cons}. 
 {\setlength{\parindent}{0cm}
\begin{thm}\label{eps_opt_theorem}
For any given $\epsilon$$>$$0$, there exists $\overline{\epsilon}$$>$$0$ such that
\begin{align}
    \min_{\pi\in\Pi^S_{R,\theta_d}(\mathcal{M})}g_{\overline{\epsilon}}(\pi)\leq \min_{\pi\in\Pi^S_{R,\theta_d}(\mathcal{M})}g_0(\pi)+\epsilon.
\end{align}
\end{thm}}
A proof of the above result and the explicit form of the constant $\overline{\epsilon}$ in terms of $\epsilon$ can be found in \cite{incentive_CDC}. The above result, together with \eqref{necessary_ineq} and \eqref{sufficient_ineq}, implies that one can find an $\epsilon$-optimal solution to the BMP-D in two steps as follows. First, solve the optimization problem on the left hand side of \eqref{sufficient_ineq} for the corresponding $\overline{\epsilon}$$>$$0$. Let $\pi^{\star}$$\in$$\Pi^S_{R,\theta_d}(\mathcal{M})$ be the optimal stationary policy which exists by Lemma \ref{stationary sufficient}. Second, provide the agent the stationary incentive sequence
\begin{align}\label{opt_incentives_dominant}
    \gamma(s,a)=\begin{cases}
    \phi_{\overline{\epsilon}}(s,a)& \ \text{if}\  \pi^{\star}(s)=a,\\
    0 &\ \text{otherwise}.
    \end{cases}
\end{align}

\subsection{Behavior modification of a dominant type: an algorithm}\label{algorithm_dominant_sec}

In the previous section, we showed that to find an $\epsilon$-optimal solution to the BMP-D, one can compute a stationary policy 
\begin{align}\label{dominant_opt_def}
\pi^{\star}\in\min_{\pi\in\Pi^S_{R,\theta_d}(\mathcal{M})}g_{\overline{\epsilon}}(\pi)
\end{align}
and provide the agent the incentive sequence given in \eqref{opt_incentives_dominant}. In this section, we present an algorithm that solves the problem in \eqref{dominant_opt_def} in time polynomial in the size of $\mathcal{M}$.



The optimization problem in \eqref{dominant_opt_def} is an instance of the so-called constrained MDP problem in which the objective is to synthesize a policy that maximizes the expected total reward while ensuring that the expected total cost is below a certain threshold \cite{Altman}. It is known that, in general,  \textit{deterministic} policies $\pi$$\in$$\Pi_{R,\theta_d}(\mathcal{M})$ are not sufficient to solve the constrained MDP problem \cite{Altman}. However, in what follows, we show that the problem in \eqref{dominant_opt_def} indeed admits an optimal solution in the set of stationary deterministic policies.

Let $r$$:$$S$$\times$$\mathcal{A}$$\rightarrow$$\mathbb{R}_{\geq 0}$ be a function such that
\begin{align}\label{reward_function}
    r(s,a):=\begin{cases}
    \sum_{s'\in B}\mathcal{P}_{s,a,s'} & \text{if} \ \ s\in S_r\\
    0 & \text{otherwise.}
    \end{cases}
\end{align}
By Theorem 10.100 in \cite{Model_checking}, for any $\pi$$\in$$\Pi(\mathcal{M})$, we have
\begin{align*}
    \mathbb{E}^{\pi}\Bigg[\sum_{t=1}^{\infty}r(S_t,\mathcal{A}_t)\Bigg]=\text{Pr}^{\pi}_{\mathcal{M}}(Reach[B]).
\end{align*}
Hence, the problem in \eqref{dominant_opt_def} can be equivalently written as
\begin{subequations}\label{main_problem2}
\begin{align}\label{probprob}
     \min_{\pi\in\Pi^S(\mathcal{M})}&\ \ \mathbb{E}^{\pi}\Bigg[\sum_{t=1}^{\infty}\phi_{\overline{\epsilon}}(S_t,\mathcal{A}_t)\Bigg]\\ \label{constraint}
     \text{subject to:}&\ \  \mathbb{E}^{\pi}\Bigg[\sum_{t=1}^{\infty}r(S_t,\mathcal{A}_t)\Bigg]=R_{\max}(\mathcal{M},B).
\end{align}
\end{subequations}

We synthesize an optimal stationary policy $\pi^{\star}$$\in$$\Pi^S(\mathcal{M})$ that solves the problem in \eqref{probprob}-\eqref{constraint} by solving two linear programs (LPs). First, we solve the LP
\begin{subequations}\label{opt_1}
\begin{flalign}\label{probobj}
&\underset{x(s,a)}{\text{minimize}}\qquad\sum_{s\in S_r}\sum_{a\in \mathcal{A}} x(s,a)\phi_{\overline{\epsilon}}(s,a)\\ \label{cons_1}
&\text{subject to:} \qquad\sum_{s\in S_r}\sum_{a\in \mathcal{A}} x(s,a)r(s,a)=R_{\max}(\mathcal{M},B)\raisetag{22pt}\\ \label{cons_2}
&\forall s\in S_r, \ \sum_{a\in \mathcal{A}(s)} x(s,a)-\sum_{s'\in S_r}\sum_{a\in \mathcal{A}(s)} \mathcal{P}_{s',a,s}x(s',a)=\alpha(s)  \\ \label{cons_last}
&\forall s\in S_r, \ a\in\mathcal{A}(s'),\  x(s,a)\geq 0  \
\end{flalign}
\end{subequations}
where $\alpha$$:$$S$$\rightarrow$$\{0,1\}$ is a function such that $\alpha(s_1)$$=$$1$ and $\alpha(s)$$=$$0$ for all $s$$\in$$S\backslash\{s_1\}$, i.e., the initial state distribution. The variable $x(s,a)$ corresponds to the expected residence time in the state-action pair $(s,a)$ \cite{Marta,Puterman}. The constraint in \eqref{cons_1} ensures that the probability of reaching the set $B$ is maximized, and the constraints in \eqref{cons_2} represent the balance between the \say{inflow} to and \say{outflow} from states.

Let $\upsilon^{\star}$ be the optimal value of the LP in \eqref{probobj}-\eqref{cons_last}. Next, we solve the  LP
\begin{subequations}\label{opt_2}
\begin{align}\label{obj22}
&\underset{x(s,a)}{\text{minimize}}\qquad\sum_{s\in S_r}\sum_{a\in \mathcal{A}} x(s,a)\\ \label{cons_2_1}
&\text{subject to:} \sum_{s\in S_r}\sum_{a\in \mathcal{A}} x(s,a)\phi_{\overline{\epsilon}}(s,a)=\upsilon^{\star}\\ 
&\qquad \qquad \qquad \eqref{cons_1}-\eqref{cons_last}. \label{const_last2}
\end{align}
\end{subequations}
Let $\{x^{\star}(s,a)$$:$$s$$\in$$S, a$$\in$$\mathcal{A}\}$ be a set of optimal decision variables for the problem in \eqref{obj22}-\eqref{const_last2}. Moreover, for a given state $s$$\in$$S$, let $\mathcal{A}^{\star}(s)$$:=$$\{a$$\in$$\mathcal{A}(s) : x^{\star}(s,a)$$>$$0 \}$ be the set of \textit{active} actions. An optimal stationary policy $\pi^{\star}$ that satisfies the condition in \eqref{dominant_opt_def} can be generated from this set as 
\begin{align}\label{arbitrary_policy_rule}
\pi^{\star}(s)=a \ \text{for an arbitrary} \ a\in\mathcal{A}^{\star}(s).
\end{align}

{\setlength{\parindent}{0cm}
\begin{prop}\label{comp_prop}
A stationary policy generated from the optimal decision variables $x^{\star}(s,a)$ using the rule in \eqref{arbitrary_policy_rule} is a solution to the problem in \eqref{probprob}-\eqref{constraint}.
\end{prop}}

A proof of Proposition \ref{comp_prop} can be found in \cite{incentive_CDC}. Intuitively, the LP in \eqref{obj22}-\eqref{const_last2} computes the minimum expected time to reach the set $B$ with probability $R_{\max}(\mathcal{M},B)$ and cost $\upsilon^{\star}$. Therefore, if $x^{\star}(s,a)$$>$$0$, by taking the action $a$$\in$$\mathcal{A}(s)$, the agent has to \say{get closer} to the set $B$ with nonzero probability. Otherwise, the minimum expected time to reach the set $B$ would be strictly decreased. Hence, by choosing an action $a$$\in$$\mathcal{A}^{\star}(s)$, the agent reaches the set $B$ with maximum probability at minimum total cost.

\section{Algorithms for Behavior Modification}\label{algorithms_section}
In this section, we present two algorithms to solve the N-BMP. In Appendix \ref{appendix_algorithm}, we explain how to modify these algorithms to solve the NS-BMP.

\subsection{Computing a globally optimal solution to the N-BMP} \label{milp_section}
We formulate a mixed-integer linear program (MILP) to compute a globally optimal solution to the N-BMP. Recall that in the N-BMP, the objective is to synthesize a \textit{stationary} incentive sequence $\gamma$$\in$$\Gamma^S(\mathcal{M})$ under which the agent reaches the target set $B$ with probability $R_{\max}(\mathcal{M},B) $ \textit{regardless of its type} and the expected total cost to the principal is minimized.

We present the MILP formulated to obtain an $\epsilon$-optimal stationary incentive sequence for the N-BMP in \eqref{MILP_begin}-\eqref{MILP_end}. In what follows, we explain its derivation in detail. Recall from Section \ref{reformulation_section} that it is sufficient for the principal to incentivize the agent only from the states $S_r$$\subseteq$$S$ in order to induce the desired behavior. Hence, we focus only on the set $S_r$ of states in the MILP formulation.

\begin{figure*}\hrule
\begin{subequations}
\begin{align}\label{MILP_begin}
    \underset{\substack{\omega, \gamma, \mathcal{V}_{\theta}, \mathcal{V}_{r,\theta}, \\ \mathcal{V}_{p,\theta}, X^{\theta}_{s,a}, \lambda_{r,\theta}  }}{\text{minimize}} &\  \omega\\
  \text{subject to:}\label{milp_1}
 &\ \forall \theta\in\Theta,  \forall s\in S_r, \forall a\in\mathcal{A}(s),  \quad \  \mathcal{Q}_{\theta}(s,a)=\mathcal{R}_{\theta}(s,a)+\gamma(s,a)\quad   \\ \label{milp_2}
  & \ \forall \theta\in\Theta,  \forall s\in S_r,\ \quad \qquad \qquad \quad \mathcal{V}_{\theta}(s)\geq\mathcal{Q}_{\theta}(s,a) \qquad\qquad \quad \qquad \qquad \qquad \quad\ \qquad \  \forall a\in\mathcal{A}(s)\\ \label{milp_3}
    & \ \forall \theta\in\Theta,  \forall s\in S_r, \ \quad \qquad  \qquad \quad \mathcal{V}_{\theta}(s)\leq\mathcal{Q}_{\theta}(s,a)+(1-X^{\theta}_{s,a})M_{\overline{\epsilon}} \ \  \qquad  \qquad \qquad \ \ \  \forall a\in\mathcal{A}(s)\\ \label{milp_4}
    & \ \forall \theta\in\Theta,  \forall s\in S_r, \forall a\in\mathcal{A}(s), \quad \   \mathcal{Q}_{\theta}(s,a)+(1-X^{\theta}_{s,a})M_{\overline{\epsilon}} \geq \mathcal{Q}_{\theta}(s,a')+\overline{\epsilon} \quad \qquad \qquad    \forall a'\in\mathcal{A}\backslash \{a\}\\ \label{milp_5}
        &\ \forall \theta\in\Theta,  \forall s\in S_r, \qquad \qquad \quad \ \ \ \sum_{a\in\mathcal{A}(s)}\mu_{r,\theta}(s,a)-\sum_{s'\in S_r}\sum_{a\in\mathcal{A}(s')}\mu_{r,\theta}(s',a)=\alpha(s),\\ \label{max_reach_milp}
    &\ \forall \theta\in\Theta,  \qquad \qquad \qquad \qquad\quad \ \  \sum_{s\in S_r}\sum_{a\in\mathcal{A}(s)}\lambda_{r,\theta}(s,a)r(s,a)=R_{\max}(\mathcal{M},B),\\  \label{milp_6}
   &\ \forall \theta\in\Theta,  \forall s\in S_r, \forall a\in\mathcal{A}(s), \quad \    \mu_{r,\theta}(s,a) \geq 0, \ \ \mu_{r,\theta}(s,a) \leq \widetilde{M}X_{s,a}^{\theta}, \ \ \mu_{r,\theta}(s,a) \leq
\lambda_{r,\theta}(s,a) \\  
    \label{milp_7}
& \ \forall \theta\in\Theta,  \forall s\in S_r, \forall a\in\mathcal{A}(s), \quad \  \mu_{r,\theta}(s,a) \geq\lambda_{r,\theta}(s,a)-(1-X_{s,a}^{\theta})\widetilde{M} \\ \label{milp_8}
& \ \forall \theta\in\Theta,  \forall s\in S_r, \forall a\in\mathcal{A}(s), \quad \ \mathcal{Q}_{p,\theta}(s,a)=\gamma(s,a)+\sum_{s'\in S}\mathcal{P}_{s,a,s'} \mathcal{V}_{p,\theta}(s'),\\ \label{milp_10}
    & \ \forall \theta\in\Theta,  \forall s\in S_r, \ \qquad \qquad \quad  \ \  \mathcal{V}_{p,\theta}(s)\geq\mathcal{Q}_{p,\theta}(s,a)-(1-X^{\theta}_{s,a})\overline{M}_{\overline{\epsilon}} \quad \quad \qquad \qquad \ \  \forall a\in\mathcal{A}(s),\\  \label{milp_11}
    & \ \forall \theta\in\Theta,  \forall s\in S_r,  \quad  \qquad \qquad  \quad  \mathcal{V}_{p,\theta}(s)\leq\mathcal{Q}_{p,\theta}(s,a)+(1-X^{\theta}_{s,a})\overline{M}_{\overline{\epsilon}}\quad  \quad \qquad \qquad\ \  \forall a\in\mathcal{A}(s), \\ \label{milp_12}
    & \  \forall \theta\in\Theta,  \qquad  \quad \  \  \ \qquad \qquad \qquad \omega \geq \mathcal{V}_{p,\theta}(s_1) \\ 
    & \ \forall \theta\in\Theta,  \forall s\in S_r, \forall a\in \mathcal{A}(s) \ \  \quad \gamma(s,a)\geq 0, \lambda_{r,\theta}(s,a)\geq 0, X^{\theta}_{s,a}\in\{0,1\}, \sum_{a\in\mathcal{A}(s)}X^{\theta}_{s,a}\geq1.\label{MILP_end}
\end{align}
\end{subequations}
\hrule
\end{figure*}

We first express the optimal behavior of an agent type $\theta$, i.e., the constraints in \eqref{policy_opt_def}-\eqref{opt_set_agent}, as a set of inequality constraints. Let $\{\gamma(s,a)$$\geq$$0 \ |\  s$$\in$$S_r , a$$\in$$\mathcal{A}\}$ be a set of variables representing the incentive offers. For each $\theta$$\in$$\Theta$, $s$$\in$$S_r$, and $a$$\in$$\mathcal{A}(s)$, let $\mathcal{V}_{\theta}(s)$$\in$$\mathbb{R}$ and $\mathcal{Q}_{\theta}(s,a)$$\in$$\mathbb{R}$ be variables such that 
\begin{subequations}
\begin{align}\label{q_val_alg}
    \mathcal{Q}_{\theta}(s,a)&=\mathcal{R}_{\theta}(s,a)+\gamma(s,a),\\ \label{v_val_alg}
    \mathcal{V}_{\theta}(s)&=\max_{a\in\mathcal{A}(s)}\mathcal{Q}_{\theta}(s,a).
\end{align}
\end{subequations}
Then, for a given state $s$$\in$$S_r$, any action $a$$\in$$\mathcal{A}(s)$ that satisfies $\mathcal{V}_{\theta}(s)$$=$$\mathcal{Q}_{\theta}(s,a)$ is optimal for the agent type $\theta$. Let $\{X^{\theta}_{s,a}$$\in$$\{0,1\} \ |\  s$$\in$$S_r , a$$\in$$\mathcal{A}, \theta$$\in$$\Theta\}$ be a set of binary variables representing whether an action $a$$\in$$\mathcal{A}(s)$ is optimal for the agent type $\theta$, i.e., $X^{\theta}_{s,a}$$=$$1$. The constraint in \eqref{v_val_alg} is not linear in $\mathcal{Q}_{\theta}(s,a)$; however, using the big-M method \cite{schrijver1998theory}, we can replace \eqref{v_val_alg} \textit{exactly} by the following set of inequalities
\begin{subequations}
\begin{flalign}\label{set_of_ineq_1}
     \mathcal{V}_{\theta}(s)&\geq\mathcal{Q}_{\theta}(s,a) \qquad\qquad \qquad\quad\ \  \forall a\in\mathcal{A}(s),&&\raisetag{23pt}\\ \label{set_of_ineq_2}
     \mathcal{V}_{\theta}(s)&\leq\mathcal{Q}_{\theta}(s,a)+(1-X^{\theta}_{s,a})M_{\overline{\epsilon}} \ \ \forall a\in\mathcal{A}(s),&& \raisetag{12pt}\\ \label{set_of_ineq_3}
     \sum_{a\in\mathcal{A}(s)}X^{\theta}_{s,a}&\geq 1&&
\end{flalign}
\end{subequations}
where $M_{\overline{\epsilon}}$ is a large constant which will be specified shortly. In particular, the constant $M_{\overline{\epsilon}}$ is chosen such that if $X^{\theta}_{s,a}$$=$$0$, then the inequality in \eqref{set_of_ineq_2} trivially holds. The constraints described above correspond to the constraints in \eqref{milp_1}-\eqref{milp_3} and \eqref{MILP_end} in the formulated MILP.

Next, we introduce a set of inequality constraints which ensure that the agent type $\theta$ has a unique optimal policy. For each $s$$\in$$S_r$ and $a$$\in$$\mathcal{A}(s)$, using the big-M method, we require
\begin{align}\label{unique_opt_cons}
    \mathcal{Q}_{\theta}(s,a)+(1-X^{\theta}_{s,a})M_{\overline{\epsilon}}\geq\mathcal{Q}_{\theta}(s,a')+\overline{\epsilon}
\end{align}
for all $a'$$\in$$\mathcal{A}(s)\backslash\{a\}$, where the constant $\overline{\epsilon}$$>$$0$ is defined in terms of the constant $\epsilon$ as shown in Theorem \ref{eps_opt_theorem}.
Note that if $X^{\theta}_{s,a}$$=$$1$, then $ \mathcal{Q}_{\theta}(s,a)$$>$$ \mathcal{Q}_{\theta}(s,a')$ for all $a'$$\in$$\mathcal{A}(s)\backslash\{a\}$. Hence, the \textit{unique} optimal policy of the agent type $\theta$ satisfies $\pi^{\star}(s)$$=$$a$. Again, the choice of $M_{\overline{\epsilon}}$ guarantees that if $X^{\theta}_{s,a}$$=$$0$, then the inequality in \eqref{unique_opt_cons} trivially hold. These constraints correspond to the constraints \eqref{milp_4} in the formulated MILP.

We now express the reachability constraint in \eqref{last_cons} as a set of equality constraints. Recall from \eqref{constraint} that the constraint in \eqref{last_cons} can be written as an expected total reward constraint with respect to the reward function $r$ defined in \eqref{reward_function}. In particular, the constraint in \eqref{cons_1}, together with the constraints in \eqref{cons_2}-\eqref{cons_last}, ensures that the agent reaches the target set $B$ with maximum probability. For each $\theta$$\in$$\Theta$, $s$$\in$$S_r$, and $a$$\in$$\mathcal{A}(s)$, let $\mu_{r,\theta}(s,a)$$\in$$\mathbb{R}$ and $\lambda_{r,\theta}(s,a)$$\in$$\mathbb{R}_{\geq 0}$ be variables such that 
\begin{subequations}
\begin{align}\label{q_val_alg_2}
    &\sum_{a\in\mathcal{A}(s)}\mu_{r,\theta}(s,a)-\sum_{s'\in S_r}\sum_{a\in\mathcal{A}(s')}\mu_{r,\theta}(s',a)=\alpha(s),\\
    &\sum_{s\in S_r}\sum_{a\in\mathcal{A}(s)}\lambda_{r,\theta}(s,a)r(s,a)=R_{\max}(\mathcal{M},B),\\ \label{lalala}
   & \mu_{r,\theta}(s,a)=\lambda_{r,\theta}(s,a)X^{\theta}_{s,a}
\end{align}
\end{subequations}
where the function $\alpha$$:$$S$$\rightarrow$$\{0,1\}$ represents the initial state distribution.
 Using the principle of optimality \cite{Puterman}, one can show that, from the initial state $s_1$$\in$$S$, the agent type $\theta$ reaches the target set $B$ with probability $R_{\max}(\mathcal{M},B)$ if and only if the agent's deterministic optimal policy expressed in terms of the variables $X_{s,a}^{\theta}$ satisfies the constraints in \eqref{q_val_alg_2}-\eqref{lalala}.

The constraint in \eqref{lalala} involves a multiplication of the binary variable $X^{\theta}_{s,a}$ with the continuous variable $\lambda_{r,\theta}(s,a)$. Using McCormick envelopes \cite{mccormick1976computability}, we can replace each constraint in \eqref{lalala} \textit{exactly} by the following inequalities 
\begin{subequations}
\begin{align}\label{relax_1}
&\mu_{r,\theta}(s,a) \geq 0 ,\\
&\mu_{r,\theta}(s,a) \leq \widetilde{M}X_{s,a}^{\theta}, \\ 
&\mu_{r,\theta}(s,a) \leq
\lambda_{r,\theta}(s,a)\\ \label{relax_2}
& \mu_{r,\theta}(s,a) \geq\lambda_{r,\theta}(s,a)-(1-X_{s,a}^{\theta})\widetilde{M}.
\end{align}
\end{subequations}
The exact value of the large constant $\widetilde{M}$ will be discussed shortly. The above inequalities ensure that if $X^{\theta}_{s,a}$$=$$0$, then $\mu_{r,\theta}(s,a)$$=$$0$, and if $X^{\theta}_{s,a}$$=$$1$, then $\mu_{r,\theta}(s,a)$$=$$\lambda_{r,\theta}(s,a)$. The constraints described above are the constraints in \eqref{milp_5}-\eqref{milp_7}.

We finally express the cost of behavior modification to the principal, i.e., the objective function in \eqref{reach_obj}, with a set of inequality constraints. For each $\theta$$\in$$\Theta$, $s$$\in$$S_r$, and $a$$\in$$\mathcal{A}(s)$, let $\mathcal{V}_{p,\theta}(s)$$\in$$\mathbb{R}$ and $\mathcal{Q}_{p,\theta}(s,a)$$\in$$\mathbb{R}$ be variables such that  
\begin{subequations}
\begin{align}\label{q_val_alg_3}
    \mathcal{Q}_{p,\theta}(s,a)&=\gamma(s,a)+\sum_{s'\in S}\mathcal{P}_{s,a,s'} \mathcal{V}_{p,\theta}(s'),\\ \label{v_val_alg_3}
    \mathcal{V}_{p,\theta}(s)&=\sum_{a\in\mathcal{A}(s)}X^{\theta}_{s,a}\mathcal{Q}_{p,\theta}(s,a).
\end{align}
\end{subequations}
Using the principle of optimality, one can show that the expected total amount of incentives paid to the agent type $\theta$ by the principal is $\mathcal{V}_{p,\theta}(s_1)$ where $s_1$$\in$$S$ is the initial state of $\mathcal{M}$. Hence, the principal's objective is to find a stationary incentive sequence that minimizes $\mathcal{V}_{p,\theta}(s_1)$ over all agent types $\theta$$\in$$\Theta$. This objective is expressed in the formulated MILP with the objective function in \eqref{MILP_begin} and the constraint in \eqref{milp_12}.

The equality constraint in \eqref{v_val_alg_3} involves bilinear terms. Using the big-M method, we can replace the constraint in \eqref{v_val_alg_3} \textit{exactly} by the following set of inequality constraints
\begin{subequations}
\begin{align}\label{set_of_ineq_4}
     \mathcal{V}_{p,\theta}(s)&\geq\mathcal{Q}_{p,\theta}(s,a)-(1-X^{\theta}_{s,a})\overline{M}_{\overline{\epsilon}} \quad \forall a\in\mathcal{A}(s),\\ \label{set_of_ineq_5}
     \mathcal{V}_{p,\theta}(s)&\leq\mathcal{Q}_{p,\theta}(s,a)+(1-X^{\theta}_{s,a})\overline{M}_{\overline{\epsilon}}\quad  \forall a\in\mathcal{A}(s), 
\end{align}
\end{subequations}
where $\overline{M}_{\overline{\epsilon}}$ is a large constant which will be precisely specified shortly. The above constraints ensure that if $X^{\theta}_{s,a}$$=$$1$, we have $\mathcal{V}_{p,\theta}(s)$$=$$\mathcal{Q}_{p,\theta}(s)$; otherwise, the corresponding inequalities trivially hold. The constraints in \eqref{milp_10}-\eqref{milp_11} in the formulated MILP then correspond to the constraints explained above.

The formulation of the MILP \eqref{MILP_begin}-\eqref{MILP_end} involves only the exact representation of the problem in \eqref{reach_obj}-\eqref{last_cons} and the constraint in \eqref{unique_opt_cons} which ensures that each agent type has a unique optimal policy. Therefore, it follows from Theorem \ref{eps_opt_theorem} that an $\epsilon$-optimal solution to the N-BMP can be obtained from the optimal solution to the formulated MILP by using its optimal decision variables $\gamma^{\star}(s,a)$ as the incentive offers.

\textit{Choosing the big-M constants:}
The formulation of the MILP \eqref{MILP_begin}-\eqref{MILP_end} requires one to specify the constants $M_{\overline{\epsilon}}$, $\widetilde{M}$, and $\overline{M}_{\overline{\epsilon}}$. They should be chosen such that the constraints involving them hold trivially when the binary variable $X^{\theta}_{s,a}$$=$$0$.

The constant $M_{\overline{\epsilon}}$ appears in \eqref{milp_3} and \eqref{milp_4}. We set 
\begin{align*}
    M_{\overline{\epsilon}}:=2\Big(\max_{\theta\in \Theta, s\in S, a\in \mathcal{A}}\mathcal{R}_{\theta}(s,a)- \min_{\theta\in \Theta, s\in S, a\in \mathcal{A}}\mathcal{R}_{\theta}(s,a)+\overline{\epsilon}\Big).
\end{align*}
The rationale behind this choice is quite simple. For a given state $s$$\in$$S_r$, the principal can make any action $a$$\in$$\mathcal{A}(s)$ uniquely optimal for any agent type $\theta$$\in$$\Theta$ by offering an incentive $\gamma(s,a)$$\leq$$\max \mathcal{R}_{\theta}(s,a)$$-$$\min \mathcal{R}_{\theta}(s,a)$$+$$\overline{\epsilon}$. Hence, we have $\mathcal{Q}_{\theta}(s,a)-\mathcal{Q}_{\theta}(s,a')$$\leq$$M_{\overline{\epsilon}}$ for any $\theta$$\in$$\Theta$, any $s$$\in$$S_r$, and any $a,a'$$\in$$\mathcal{A}(s)$. Consequently, for the above choice of $M_{\overline{\epsilon}}$, the constraints in \eqref{milp_3} and \eqref{milp_4} hold trivially when $X^{\theta}_{s,a}$$=$$0$.

The constant $\widetilde{M}$ appears in \eqref{milp_6}-\eqref{milp_7}. These constraints put upper bounds on the expected residence times in state-action pairs. For MDPs with deterministic transition functions, we simply set $\widetilde{M}$$=$$1$. The reason is that, under a deterministic policy that reaches the set $B$ with maximum probability, the agent cannot visit a state $s$$\in$$S_r$ twice. On the other hand, for MDPs with stochastic transition functions, it is, in general, NP-hard to compute the maximum expected residence times in states subject to the reachability constraint in \eqref{last_cons}. Therefore, in that case, we set $\widetilde{M}$$=$$k\lvert S\rvert$ for some large $k$$\in$$\mathbb{N}$.

The constant $\overline{M}_{\overline{\epsilon}}$ appears in \eqref{set_of_ineq_4}-\eqref{set_of_ineq_5}. These constraints put bounds on the total expected cost $\mathcal{V}_{p,\theta}(s)$ to the principal. Using the methods described in Section \ref{algorithm_dominant_sec}, one can compute a policy that minimizes the expected total cost with respect to the cost function $C$, defined in \eqref{cost_cost_cost}, subject to the reachability constraint in \eqref{last_cons}. The optimal value of this optimization problem corresponds to a total incentive amount that is \textit{sufficient} to lead the agent to the target set regardless of its type. Hence, it provides an upper bound on the variables $\mathcal{Q}_{p,\theta}(s,a)$. We set the constant $\overline{M}_{\overline{\epsilon}}$ to the optimal value of the above mentioned optimization problem. As a result, the constraints in  \eqref{set_of_ineq_4}-\eqref{set_of_ineq_5} hold trivially when $X^{\theta}_{s,a}$$=$$0$.

\subsection{Computing a locally optimal solution to the N-BMP}\label{local_opt_sec}
The MILP formulation described in the previous section involves $\lvert S \rvert\lvert \mathcal{A} \rvert\lvert \Theta\rvert$ binary variables. In the worst case, computing its optimal solution takes exponential time in the size of $\mathcal{M}$ and $\Theta$. Here, we present a more practical algorithm which computes a locally optimal solution to the N-BMP. Specifically, we formulate the N-BMP as a nonlinear optimization problem (NLP) with bilinear constraints by slightly modifying the MILP \eqref{MILP_begin}-\eqref{MILP_end}. We then use the so-called convex-concave procedure (CCP) to solve the derived NLP. 

To formulate the N-BMP as an NLP with bilinar constraints instead of an MILP, we express the policy of the agent type $\theta$ using the set of continuous variables 
\begin{align*}
    \Bigg\{0\leq\nu^{\theta}(s,a)\leq1 \Bigg|  \sum_{a\in\mathcal{A}(s)}\nu^{\theta}(s,a)=1, s\in S_r, a\in \mathcal{A}, \theta\in\Theta\Bigg\}
\end{align*}
instead of the binary variables $X^{\theta}_{s,a}$$\in$$\{0,1\}$. In what follows, we explain the modifications we make to the MILP \eqref{MILP_begin}-\eqref{MILP_end} to obtain an NLP with bilinear constraints.

To express the optimal behavior of an agent type $\theta$, we keep the constraints \eqref{milp_1}-\eqref{milp_2}, but \textit{replace} the constraint in \eqref{milp_3} in the MILP with the following constraint
\begin{align}
    \mathcal{V}_{\theta}(s)&\leq\sum_{a\in\mathcal{A}(s)}\nu^{\theta}(s,a)\mathcal{Q}_{\theta}(s,a).\label{first_bilinear}
\end{align}
The above inequality, together with \eqref{milp_1}-\eqref{milp_2}, ensures that if $\nu^{\theta}(s,a)$$>$$0$, the corresponding action $a$$\in$$\mathcal{A}(s)$ is optimal for the agent type $\theta$. Note that the condition $\nu^{\theta}(s,a)$$>$$0$ is \textit{not} necessary for the action $a$$\in$$\mathcal{A}(s)$ to be optimal. Next, we introduce a constraint which ensures that the agent type $\theta$ has a unique deterministic optimal policy; therefore, the condition $\nu^{\theta}(s,a)$$>$$0$ becomes a necessary and sufficient condition for the optimality of the action $a$$\in$$\mathcal{A}(s)$ for the agent type $\theta$.

Now, similar to \eqref{unique_opt_cons}, we introduce a set of inequality constraints which ensure that the agent type $\theta$ has a unique optimal policy. In particular, we \textit{replace} the constraints in \eqref{milp_4} in the MILP with the following set of constraints
\begin{align}\label{unique_opt_cons_2}
    \nu^{\theta}(s,a)\Big[\mathcal{Q}_{\theta}(s,a)-\Big(\mathcal{Q}_{\theta}(s,a')+\overline{\epsilon}\Big)\Big]\geq 0
\end{align}
for all $a'$$\in$$\mathcal{A}(s)\backslash\{a\}$. The above constraints, together with the constraints in \eqref{milp_1}-\eqref{milp_2} and \eqref{first_bilinear}, ensure that $\nu^{\theta}(s,a)$$=$$1$ for the uniquely optimal action $a$$\in$$\mathcal{A}(s)$. 

In the MILP formulation, we express the reachability constraint in \eqref{last_cons} as the set of constraints in \eqref{milp_5}-\eqref{milp_7}. In particular, the constraints in \eqref{relax_1}-\eqref{relax_2} represent the McCormick envelope corresponding to the constraint in \eqref{lalala}. In the NLP formulation, we keep the constraints in \eqref{milp_5}-\eqref{max_reach_milp}, but instead of using the McCormick envelopes, we \textit{replace} the constraint in \eqref{lalala} simply with the following bilinear constraint
\begin{align}\label{max_reac_111}
\mu_{r,\theta}(s,a)=\lambda_{r,\theta}(s,a)\nu^{\theta}(s,a).
\end{align}
The above constraint, together with the constraints in \eqref{milp_5}-\eqref{max_reach_milp}, ensures that the agent type $\theta$$\in$$\Theta$ reaches the target set $B$ with probability $R_{\max}(\mathcal{M},B)$ under its optimal policy expressed in terms of the variables $\nu^{\theta}(s,a)$.  

Finally, we express the expected total cost of behavior modification with a set of inequality constraints. Specifically, we keep the constraints in \eqref{milp_8} and \eqref{milp_12} in the MILP, but \textit{replace} the constraints in \eqref{milp_10}-\eqref{milp_11} with the constraint
\begin{align}\label{last_bilinear}
    \mathcal{V}_{p,\theta}(s)\geq\sum_{a\in\mathcal{A}(s)}\nu^{\theta}(s,a)\mathcal{Q}_{p,\theta}(s,a).
\end{align}
It is important to note that we removed the constraint that puts an upper bound on the value of $\mathcal{V}_{p,\theta}(s)$. 
As a result, unlike the MILP \eqref{MILP_begin}-\eqref{MILP_end}, in the formulated optimization problem, the value of $ \mathcal{V}_{p,\theta}(s_1)$ \textit{is, in general, not equal to} the cost of behavior modification for the corresponding agent type $\theta$. However, it \textit{is indeed equal to} the cost of behavior modification for the agent type that maximizes the cost to the principal, i.e., $\theta$$\in$$\Theta$ that satisfies $w$$=$$\mathcal{V}_{p,\theta}(s_1)$. 

The formulation of the NLP described above involves the exact representation of the constraints in \eqref{policy_opt_def}-\eqref{last_cons} and an additional constraint in \eqref{unique_opt_cons_2} which ensures that each agent type has a unique optimal policy. Moreover, it exactly represents the worst-case expected total cost incurred by the principal. Hence, it follows from Theorem \ref{eps_opt_theorem} that an $\epsilon$-optimal solution to the N-BMP can be obtained from the optimal solution to the formulated NLP by using its optimal decision variables $\gamma^{\star}(s,a)$ as the incentive offers to the agent.

\textit{CCP to compute a local optimum:} We employ the Penalty CCP (P-CCP) algorithm \cite{lipp2016variations}, which is a variation of the basic CCP algorithm \cite{yuille2003concave}, to compute a locally optimal solution to the formulated NLP. We now briefly cover the main working principles of the P-CCP algorithm and explain how to utilize it for the purposes of this paper. We refer the interested reader to \cite{lipp2016variations} for further details on the P-CCP algorithm.

Suppose we are given an optimization problem of the form 
\begin{subequations}
\begin{align}\label{CCP_opt_1}
    \underset{x\in\mathbb{R}^n}{\text{minimize}}& \ Z_0(x)\\ \label{CCP_opt_2}
    \text{subject to:}& \ Z_i(x)-Y_i(x)\leq 0, \ \ i=1,2,\ldots,m
\end{align}
\end{subequations}
where $Z_i$$:$$\mathbb{R}^n$$\rightarrow$$\mathbb{R}$ and $Y_i$$:$$\mathbb{R}^n$$\rightarrow$$\mathbb{R}$ are \textit{convex} functions for each $i$$\in$$\{0,1\ldots,m\}$. The above optimization problem is, in general, not convex due to the constraints in \eqref{CCP_opt_2} \cite{boyd2004convex}. The P-CCP algorithm is a heuristic method to compute a locally optimal solution to the problems of the form \eqref{CCP_opt_1}-\eqref{CCP_opt_2} by iteratively approximating the constraints in \eqref{CCP_opt_2}. 

Let $\zeta$ be a constant such that $\zeta$$>$$1$, $\tau_0$ and $\tau_{\max}$ be positive constants, and $\tau_k$ be a constant that is recursively defined as $\tau_k$$:=$$\min\{\zeta\tau_{k-1},\tau_{\max}\}$ for $k$$\in$$\mathbb{N}$. Note that $\tau_i$$\geq$$\tau_j$ for all $i$$\geq$$j$.

Starting from an arbitrary initial point $x_0$$\in$$\mathbb{R}^n$, at the $k$-th iteration, the P-CCP algorithm computes a globally optimal solution $x_{k+1}$ of the \textit{convex optimization problem}
\begin{subequations}
\begin{align}\label{CCP_opt_3}
    \underset{x\in\mathbb{R}^n}{\text{minimize}}& \ Z_0(x)+\tau_k\sum_{i=1}^m \overline{s}_i^k\\ \label{CCP_opt_4}
    \text{subject to:}& \ Z_i(x)-\overline{Y}_i(x;x_k)\leq \overline{s}_i^k, \ \ i=1,2,\ldots,m\\ \label{CCP_opt_5}
    & \  \overline{s}_i^k\geq 0, \ \ i=1,2,\ldots,m
\end{align}
\end{subequations}
where $\overline{Y}_i(x;x_k)$$:=$$Y_i(x_k)$$+$$\nabla Y_i(x_k)^T(x$$-$$x_k)$ is the first-order approximation of the function $Y_i$ at $x_k$, and $\nabla Y_i$ denotes the gradient of $Y_i$. In the above problem, the variables $\overline{s}_i^k$ is referred to as the \textit{slack variables}. This is because, when the optimal solution to the problem in \eqref{CCP_opt_3}-\eqref{CCP_opt_5} satisfies $\overline{s}_i^k$$=$$0$ for all $i$, then the optimal variables $x_{k+1}$$\in$$\mathbb{R}^n$ constitute a feasible solution for the original problem in \eqref{CCP_opt_1}-\eqref{CCP_opt_2}.

The P-CCP algorithm terminates when 
\begin{align*}
    \Big(Z_0(x_k)+\tau_k\sum_{i=1}^m \overline{s}_i^k\Big)-\Big(Z_0(x_{k+1})+\tau_k\sum_{i=1}^m \overline{s}_i^{k+1}\Big)\leq \overline{\delta}
\end{align*}
for some small $\overline{\delta}$$>$$0$, and either $x_k$ is feasible, i.e., 
\begin{align}\label{feasibility_condition_1}
    \sum_{i=1}^m \overline{s}_i^{k+1}\leq \overline{\delta}_{\text{violation}}\approx 0,
\end{align}
or $\tau_k$$=$$\tau_{\max}$. The P-CCP algorithm is guaranteed to terminate \cite{lipp2016variations}. Moreover, if the condition in \eqref{feasibility_condition_1} is satisfied upon termination, the output of the algorithm constitutes a locally optimal solution to the problem in \eqref{CCP_opt_1}-\eqref{CCP_opt_2}. 

The NLP formulated in the previous section consists of a convex objective function, i.e., $Z_0(x)$$=$$w$, a number of linear constraints, and the bilinear constraints in \eqref{first_bilinear}-\eqref{last_bilinear}. To employ the P-CCP algorithm for obtaining a locally optimal solution to the formulated NLP, we express the above mentioned bilinear constraints in the form of \eqref{CCP_opt_2} by following the convexification technique described in \cite{cubuktepe2018synthesis}. Due to space restrictions, we refer the interested reader to Section 5 in \cite{cubuktepe2018synthesis} for further details on the convexification technique.

We can now compute a locally optimal solution to the NLP formulated in the previous section as follows. Set the parameters $\zeta$, $\tau_0$, $\tau_{\max}$, $\overline{\delta}$, and $\overline{\delta}_{\text{violation}}$. Initialize the P-CCP algorithm by setting all variables to some initial values. Run the P-CCP algorithm by convexifying each bilinear term in the NLP at each iteration. Upon termination, verify whether the condition in \eqref{feasibility_condition_1} is satisfied. If it is satisfied, use the optimal decision variables $\gamma^{\star}(s,a)$ as the incentive offers to the agent.
\section{Numerical Examples}
We illustrate the application of the presented algorithms with two examples on discount planning and motion planning.\footnote{ Due to space restrictions, we provide only a brief description of the reward functions in numerical examples. We refer the interested readers to \url{https://github.com/yagizsavas/Sequential-Incentive-Design} for a complete description of reward functions and an implementation of the algorithms.} We run the computations on a 3.1 GHz desktop with 32 GB RAM and employ the GUROBI solver \cite{gurobi} for optimization.

\subsection{Discount planning to encourage purchases}
We consider a retailer (principal) that aims to sell $n$ products to a customer (agent) at the end of $n$ interactions. It is assumed that the agent purchases a single product per interaction, and the agent's willingness to purchase new products depends on the ones purchased in the past. The principal does not know how the agent associates different products with each other and aims to maximize its total profit by synthesizing a sequence of discount offers (incentives) by convincing the agent to purchase all the products with a minimum total discount. 

Let $Q$$=$$\{1,2,\ldots,n\}$ be the set of products that the principal aims to sell. We construct an MDP to express the agent's behavior as follows. Each state in the MDP corresponds to a set of products that are already purchased by the agent, i.e., the set $S$ of states is the power set of $Q$. The initial state of the MDP is the empty set. In a given state $s$$\in$$S$ such that $s$$\subseteq$$Q$, the agent has two choices: not to make any purchases and stay in the same state by taking the action $a_0$ or to purchase a product $i$$\in$$Q\backslash s$ by taking the action $a_i$ as a result of which it transitions to the state $s\cup\{i\}$$\in$$S$. All transitions in this model are deterministic, i.e., $\mathcal{P}_{s,a,s'}$$\in$$\{0,1\}$. 

We choose $n$$=$$4$ and consider 3 agent types, i.e., $\lvert \Theta\rvert$$=$$3$. All agent types prefer not to purchase any products unless the principal offers some discounts, i.e., $\mathcal{R}^{\max}_{\theta}(s)$$=$$\mathcal{R}_{\theta}(s,a_0)$$=$$0$ for all $s$$\in$$S$ and $\theta$$\in$$\Theta$. The intrinsic motivation of each agent type is summarized in Table \ref{specs_numeric}. The product grouping in Table \ref{specs_numeric} expresses how an agent type associates the products with each other. For example, once the agent type $\theta_1$ purchases product $1$, it becomes more willing to purchase product $2$ than products $3$ or $4$, i.e., $\mathcal{R}_{\theta_1}(\{1\},a_2)$$=$$-1$ and $\mathcal{R}_{\theta_1}(\{1\},a_3)$$=$$\mathcal{R}_{\theta_1}(\{1\},a_4)$$=$$-2$. Once a group of products is purchased by the agent, it becomes harder for the principal to sell another group of products. The product importance in Table \ref{specs_numeric} expresses the agent's willingness to purchase a group of products once the other group of products is already purchased. For example, it is more profitable for the principal to sell first the group $G_1$$=$$\{1,2\}$ of products to the type $\theta_1$ and then the group $G_2$$=$$\{3,4\}$ of products since $\theta_1$ is already more willing to purchase the group $G_2$, i.e., $\mathcal{R}_{\theta_1}(\{1,2\},a_3)$$=$$\mathcal{R}_{\theta_1}(\{1,2\},a_4)$$=$$-2$ and $\mathcal{R}_{\theta_1}(\{3,4\},a_1)$$=$$\mathcal{R}_{\theta_1}(\{3,4\},a_2)$$=$$-3$. 

\begin{table}[t!]\centering
\caption{Customer preferences in discount planning.}
\begin{tabular}{ |c || c|| c| } 
 \hline
 {{\textbf{Customer Type}}} & \textbf{Product grouping} & \textbf{Product importance} \\ 
 \hline
   \hline 
  $\theta_1$ & $G_1$$=$$\{1,2\}$, $G_2$$=$$\{3,4\}$ & $G_2$$>$$G_1$ \\ 
 \hline
 $\theta_2$ & $G_1$$=$$\{1,3\}$, $G_2$$=$$\{2,4\}$ & $G_2$$>$$G_1$  \\ 
  \hline
  $\theta_3$ & $G_1$$=$$\{1,4\}$, $G_2$$=$$\{2,3\}$ & $G_2$$=$$G_1$ \\ 
 \hline
\end{tabular}
\label{specs_numeric}
\end{table}

We synthesized a sequence of discount offers for the agent by solving the corresponding MILP formulation which has 205 continuous and 112 binary variables. The computation took 0.4 seconds. A part of the synthesized discount sequence is demonstrated in Table \ref{discount_results}. During the first interaction, the principal offers a discount only for product 1. The offer is large enough to ensure that the agent purchases product 1 regardless of its type, i.e., $\gamma(\{\},a_1)$$=$$1$$+$$\epsilon$. Then, during the second interaction, the principal offers discounts $\gamma(\{1\},a_2)$$=$$\gamma(\{1\},a_3)$$=$$\gamma(\{1\},a_4)$$=$$1$$+$$\epsilon$ for all remaining products. As a result, depending on its type, the agent purchases one of the discounted products. For example, the agent type $\theta_1$ purchases product 2, whereas the agent type $\theta_2$ purchases product 3. Finally, the principal sequentially offers discounts for remaining two products.

The synthesized discount sequence utilizes the principal's knowledge on the possible agent types. Specifically, the principal knows that purchasing product 1 is not a priority for any of the agent types (see Table \ref{specs_numeric}). Therefore, the principal is aware that it will be harder to sell product 1 once the agent spends money on the other products. Hence, by offering a discount only for product 1 during the first interaction, the principal makes sure that the agent purchases the least important product when it still has the money. Moreover, during the second interaction, the principal ensures that each agent type purchases the second least important product for itself by offering discounts for products 2, 3, and 4 at the same. Since the group of products that are not purchased yet are more important than the group of products that are already purchased, the principal then sells the remaining products by offering small discount amounts.

\subsection{Incentives for motion planning }\label{motion_plan_sec}
\begin{figure}[t]
\begin{subfigure}[b]{0.3\linewidth}
\scalebox{0.23}{
\includegraphics[]{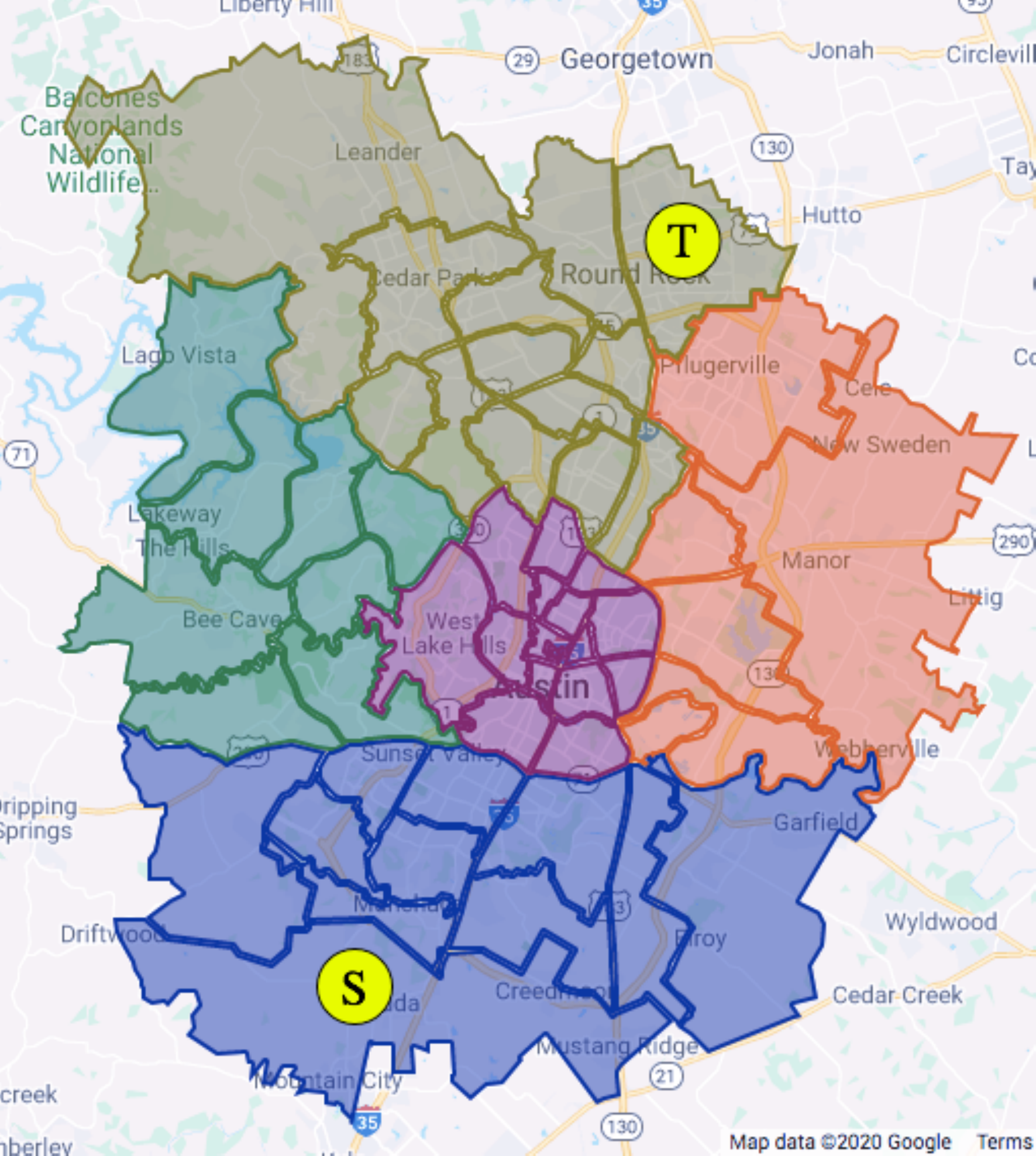}}
\end{subfigure}
\hspace{0.2\linewidth}
\begin{subfigure}[b]{0.3\linewidth}
\scalebox{0.23}{
\includegraphics[]{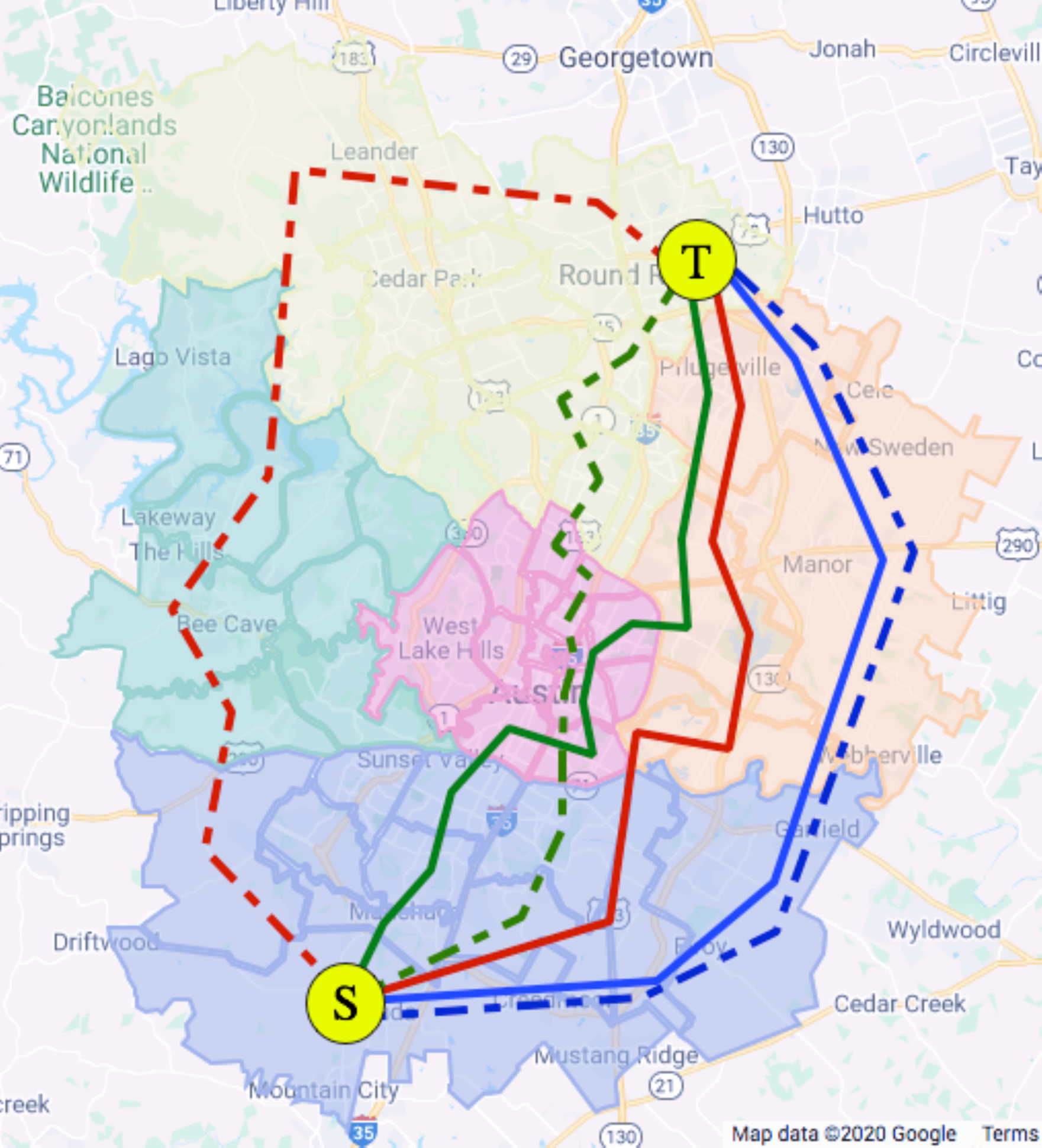}}
\end{subfigure}
\caption{Motion planning with incentive offers. The states with labels $S$ and $T$ are the initial and target states, respectively. Green, blue, and red lines are the trajectories followed by the agent types $\theta_1$, $\theta_2$, and $\theta_3$, respectively. When the agent type is unknown (known) to the principal, the solid (dashed) trajectories are followed under the synthesized incentive sequence. }
\label{grid_graph}
\end{figure}

We consider a ridesharing company (principal) that aims to convince a driver (agent) to be present at a desired target region during the rush hour by offering monetary incentives. The company is assumed to operate in Austin, TX, USA which is divided into regions according to zip codes as shown in Fig. \ref{grid_graph} (left) in which the areas, e.g., downtown, north etc., are indicated with different colors for illustrative purposes. The environment is modeled as an MDP with 54 states each of which corresponds to a region. In state $s_i$$\in$$S$, the agent has two choices: to stay in the same state $s_i$ under the action $a_i$ or to move to a ``neighboring" state $s_j$, i.e., a region that share a border with the agent's current region, under the action $a_j$. 
\begin{table}[t!]\centering
\caption{Discount offers that maximize the retailer's profit while ensuring that the customer purchases all products regardless of its type. The check marks indicate the discounted products as a function of the products that are already purchased. For example, the retailer discounts only product 1 if the customer has no previous purchases and discounts product 4 if the set $\{1,2\}$ of products has already been purchased. }
\scalebox{0.9}{
\begin{tabular}{ |c | c| c| c| c| c|} 
\hline
\diagbox{\textbf{Product}}{\textbf{Purchased}} & $\{\}$ & $\{1\}$ & $\{1,2\}$ & $\{1,3\}$ & $\{1,4\}$  \\ \hline
1 & \checkmark & - & - & - & -  \\ \hline
2 & - & \checkmark & - & - & -   \\ \hline
3 & - & \checkmark &  - & - & \checkmark \\ \hline
4 & - & \checkmark & \checkmark & \checkmark & - \\ \hline
\end{tabular}}
\label{discount_results}
\end{table}

We consider three agent types, i.e., $\lvert \Theta \rvert$$=$$3$. In the absence of incentives, staying in the same state is optimal for all agent types in all states, i.e., $\mathcal{R}^{\max}_{\theta}(s_i)$$=$$\mathcal{R}_{\theta}(s_i,a_i)$$=$$0$ for all $s_i$$\in$$S$ and $\theta$$\in$$\Theta$. Hence, the principal cannot distinguish the true agent type by passively observing the agent and, to induce the desired behavior, the principal has to offer incentives in all states. The first agent type $\theta_1$ associates each state pair $(s_i,s_j)$ with their corresponding distance $d_{i,j}$. To move to a different region, it demands the principal to offer incentives that is proportional to $d_{i,j}$, i.e., $\mathcal{R}_{\theta_1}(s_i,a_j)$$=$$-d_{i,j}$. The second agent type $\theta_2$ associates each state $s_i$ with a congestion index $t_i$$\in$$\{1,2,\ldots,10\}$, e.g., the states in $downtown$ area has the highest congestion indices. To move to a region $s_j$, it demands the principal to offer incentives that is proportional to $t_{j}$, i.e., $\mathcal{R}_{\theta_2}(s_i,a_j)$$=$$-2t_{j}$. Finally, the third agent type $\theta_3$ takes both the distance and congestion index into account and has the reward function $\mathcal{R}_{\theta_3}(s_i,a_j)$$=$$-(0.8d_{i,j}+0.2t_{j})$.

We first considered the case of \textit{known} agent types and synthesized optimal incentive sequences for \textit{each} agent type using the corresponding LP formulation. Trajectories followed by the agents under the synthesized incentive sequences are shown in Fig. \ref{grid_graph} (right) with dashed lines, e.g., the type $\theta_1$ is incentivized to follow the shortest trajectory to the target. We then considered the case of $\textit{unknown}$ agent types and computed a globally optimal solution to the corresponding N-BMP instance using the MILP formulation. The MILP had 1407 continuous and 937 binary variables, and the computation exceeded the memory limit after 6 hours. Finally, we computed a locally optimal solution to the corresponding N-BMP instance through the NLP formulation which has 64801 continuous variables. The CCP converged in 193 iterations, and the computation took 1543 seconds in total. Optimal trajectories followed by each agent type under the synthesized incentive sequence are shown in Fig. \ref{grid_graph} (right) with solid lines. We measured the suboptimality of the synthesized incentive sequence using the lower bound in \eqref{weak_duality}; the total cost of the synthesized incentive sequence to the principal is at most 1.52 times the total cost of the globally optimal one.

As seen in Fig. \ref{grid_graph}, under the synthesized incentive sequence, the agent reaches the target state desired by the principal regardless of its type. Moreover, the synthesized incentive sequence utilizes the principal's information on possible agent types. In particular, under the offered incentive sequences, the agent types $\theta_1$ and $\theta_2$ follow trajectories that resemble the ones that would be followed if the true agent type was known to the principal. Moreover, even though the optimal trajectory of the type $\theta_3$ changes significantly due to the principal's incomplete information on the true agent type, the type $\theta_3$ still avoids the downtown area under the offered incentive sequences.

\newcommand{\StaticObstacle}[2]{ \fill[red] (#1+0.1,#2+0.1) rectangle (#1+0.9,#2+0.9);}
\newcommand{\initialstate}[2]{ \fill[brown] (#1+0.15,#2+0.15) rectangle (#1+0.85,#2+0.85);}
\newcommand{\goalstate}[2]{ \fill[green] (#1+0.15,#2+0.15) rectangle (#1+0.85,#2+0.85);}
\section{Conclusions and Future Directions}
We considered a principal that offers a sequence of incentives to an agent with unknown intrinsic motivation to align the agent's behavior with a desired objective. We showed that the behavior modification problem (BMP), the problem of synthesizing an incentive sequence that induces the desired behavior at minimum total cost is, in general, computationally intractable. We presented a sufficient condition under which the BMP can be solved in time polynomial in the related parameters. Finally, we developed two algorithms to synthesize a stationary incentive sequence that induces the desired behavior while minimizing the total cost either globally or locally.

In the BMP, we assume that the agent's intrinsic motivation can be expressed as a reward function that belongs to a finite set of possible reward functions. The performance of any algorithm that solves the BMP necessarily depends on the size of this finite set. Hence, it is of interest to develop methods that extract reward representations from data by taking the computational complexity of the resulting BMP into account. Another possible future direction is to express the agent's unknown intrinsic motivation as a generic function with certain structural properties, e.g., continuity, monotonicity, concavity etc., and investigate the effects of different properties on the complexity of the resulting BMP.
  \bibliographystyle{IEEEtran}
\bibliography{main.bib}
\appendices
\section{}\label{Appendix_remarks}
In this appendix, we briefly explain how to modify the methods developed in this paper in order to solve the BMP when the agent has a finite decision horizon $N$$\in$$\mathbb{N}$ and when the set $B$$\subseteq$$S$ of target states consists of non-absorbing states. 

Consider an agent with a finite decision horizon $N$$\in$$\mathbb{N}$ whose objective is to maximize the expected total reward it collects at the end of \textit{every} $N$ stages \cite{incentive_CDC}. Such an agent's optimal policy is a sequence $(D_1,D_2,\ldots)$ of decisions where, for $k$$\in$$\mathbb{N}$, $D_k$$:=$$(d_{(k-1)N+1},d_{(k-1)N+2},\ldots,d_{kN})$ satisfies
\begin{align*}
 D_k \in \arg\max_{\pi\in\Pi(\mathcal{M})} \mathbb{E}^{\pi}\Bigg[\sum_{t=(k-1)N+1}^{kN}\mathcal{R}_{\theta^{\star}}(s_t,a)+\delta_t(I_t,a)\Bigg].
\end{align*}

We now show that the behavior modification of an agent with a finite decision horizon $N$ is computationally not easier than the behavior modification of a myopic agent. For any MDP $\mathcal{M}$ and $N$$\in$$\mathbb{N}$, consider an \textit{expanded MDP} with the set $\overline{S}$$:=$$S$$\times$$[N]$ of states and the initial state $(s_1,1)$. Let the transition function $\overline{\mathcal{P}}$$:$$\overline{S}$$\times$$\mathcal{A}$$\times$$\overline{S}$$\rightarrow$$[0,1]$ be such that $\overline{\mathcal{P}}((s,t),a,(s,t$$+$$1))$$=$$1$ for $t$$\in$$[N]\backslash\{N\}$ and $\overline{\mathcal{P}}((s,N),a,(s',1))$$=$$\mathcal{P}_{s,a,s'}$, and the reward function $\overline{\mathcal{R}}_{\theta^{\star}}$$:$$\overline{S}$$\times$$ \mathcal{A}$$\rightarrow$$\mathbb{R}$ be such that $\overline{\mathcal{R}}_{\theta^{\star}}((s,N),a)$$:=$$\mathcal{R}_{\theta^{\star}}(s,a)$ and $\overline{\mathcal{R}}_{\theta^{\star}}((s,t),a)$$:=$$0$ otherwise. It can be shown that, on the expanded MDP, the behavior modification of an agent with a decision horizon $N$ is equivalent to the behavior modification of a myopic agent. Since the expanded MDP is constructed in time polynomial in the size of $\mathcal{M}$ and $N$, the result follows.

One can solve the BMP when the agent has a decision horizon $N$$\in$$\mathbb{N}$ as follows. For a given BMP instance, first construct an MDP with the set $\overline{S}$$:=$$S$$\times$$[N]$ of states, the initial state $(s_1,1)$, and the transition function $\overline{\mathcal{P}}$$:$$\overline{S}$$\times$$\mathcal{A}$$\times$$\overline{S}$$\rightarrow$$[0,1]$ such that $\overline{\mathcal{P}}((s,t),a,(s',t$$+$$1))$$=$$\mathcal{P}_{s,a,s'}$ for $t$$\in$$[N]\backslash\{N\}$ and $\overline{\mathcal{P}}((s,N),a,(s',1))$$=$$\mathcal{P}_{s,a,s'}$. Note that, in the above construction, the agent returns to a state $(s,1)$ where $s$$\in$$S$ after every $N$ stages. Second, on the constructed MDP, synthesize a sequence of incentive offers that solve the BMP. Finally, in every $N$ stages, provide the agent the next $N$ incentive offers. 

One can solve the BMP when the set $B$ of target states consists of non-absorbing states as follows. For a given BMP instance, construct an MDP with the set $\overline{S}$$:=$$S$$\times$$\{1,2\}$ of states, the initial state $(s_1,1)$, and the transition function $\overline{\mathcal{P}}$$:$$\overline{S}$$\times$$\mathcal{A}$$\times$$\overline{S}$$\rightarrow$$[0,1]$ such that $\overline{\mathcal{P}}((s,1),a,(s',1))$$=$$\mathcal{P}_{s,a,s'}$ for all $s$$\in$$S\backslash B$, $\overline{\mathcal{P}}((s,1),a,(s',2))$$=$$\mathcal{P}_{s,a,s'}$ for all $s$$\in$$B$, and $\overline{\mathcal{P}}((s,2),a,(s',2))$$=$$\mathcal{P}_{s,a,s'}$ for all $s$$\in$$S$. On the constructed MDP, the agent transitions to a state $(s,2)$$\in$$S$$\times$$\{2\}$ if and only if it reaches the target set. Make all states $S$$\times$$\{2\}$ absorbing and replace the target set $B$ with the set $S$$\times$$\{2\}$. Finally, on the constructed MDP, synthesize a sequence of incentive offers that solve the BMP.
\section{}\label{appendix_proofs}
In this appendix, we provide proofs for all results presented in this paper. 

{\setlength{\parindent}{0cm}
\begin{definition}
\textbf{(QSAT)} \cite{papadimitriou1987complexity} Let $F(x_1,x_2,\ldots,x_n)$ be a Boolean formula in conjunctive normal form with three literals per clause, and $\exists x_1 \forall x_2 \exists x_3 \ldots \forall x_n F(x_1,x_2,\ldots,x_n)$ be a quantified Boolean formula (QBF). Decide whether the given QBF is true.
\end{definition}}

\noindent \textbf{Proof of Theorem \ref{PSPACE_thm}:} The proof is by a reduction from QSAT. We are given an arbitrary QBF with $n$ variables and $m$ clauses $C_1,C_2,\ldots,C_m$. Without loss of generality, we assume that $n$ is an even number. To prove the claim, we first construct an MDP $\mathcal{M}$, a target set $B$$\subseteq$$S$, a type set $\Theta$, and a reward function $\mathcal{R}_{\theta}$ for each $\theta$$\in$$\Theta$. We then show that, on the constructed model, the total cost of an incentive sequence that solves the BMP is $n$ or less if and only if the QBF is true.

\begin{figure*}[t!]\centering
\scalebox{0.8}{
\begin{tikzpicture}[->, >=stealth', auto, semithick, node distance=2cm]

    \tikzstyle{every state}=[fill=white,draw=black,thick,text=black,scale=0.8]

    \node[state,initial, initial text=, fill=green] (A_1) {$A_1'$};
    \node[state,fill=green] (A_11) [ below =25mm of A_1]  {$A_1$};
    
    \node[state,fill=red] (T_1) [above right =5mm and 10 mm of A_1]  {$T_1'$};
    \node[state,fill=red] (F_1) [below right =5mm and 10 mm of A_1]  {$F_1'$};

    \node[state,fill=red] (T_11) [above right =5mm and 10 mm of A_11]  {$T_1$};
    \node[state,fill=red] (F_11) [below right =5mm and 10 mm of A_11]  {$F_1$};
    
    \node[state,fill=orange] (A_2) [right=25 mm of A_1] {$A_2'$};
    \node[state,fill=orange] (A_22) [ below =25mm of A_2]  {$A_2$};
    
    \node[state,fill=red] (T_2) [above right =5mm and 10 mm of A_2]  {$T_2'$};
    \node[state,fill=red] (F_2) [below right =5mm and 10 mm of A_2]  {$F_2'$};

    \node[state,fill=red] (T_22) [above right =5mm and 10 mm of A_22]  {$T_2$};
    \node[state,fill=red] (F_22) [below right =5mm and 10 mm of A_22]  {$F_2$};
    
    
       \node[state,fill=green] (A_3) [right=25 mm of A_2] {$A_3'$};
    \node[state,fill=green] (A_33) [ below =25mm of A_3]  {$A_3$};
    
    \node[state,fill=red] (T_3) [above right =5mm and 10 mm of A_3]  {$T_3'$};
    \node[state,fill=red] (F_3) [below right =5mm and 10 mm of A_3]  {$F_3'$};

    \node[state,fill=red] (T_33) [above right =5mm and 10 mm of A_33]  {$T_3$};
    \node[state,fill=red] (F_33) [below right =5mm and 10 mm of A_33]  {$F_3$};
 
     \node[state,fill=orange] (A_4) [right=25 mm of A_3] {$A_4'$};
    \node[state,fill=orange] (A_44) [ below =25mm of A_4]  {$A_4$};
    
    
         \node[state,fill=orange] (A_n) [right=10 mm of A_4] {$A_n'$};
    \node[state,fill=orange] (A_n1) [ below =25mm of A_n]  {$A_n$};
    
    \node[state,fill=red] (T_n) [above right =5mm and 10 mm of A_n]  {$T_n'$};
    \node[state,fill=red] (F_n) [below right =5mm and 10 mm of A_n]  {$F_n'$};

    \node[state,fill=red] (T_n1) [above right =5mm and 10 mm of A_n1]  {$T_n$};
    \node[state,fill=red] (F_n1) [below right =5mm and 10 mm of A_n1]  {$F_n$};
    
        \node (T_n111) [below right =10mm and 2 mm of A_4,scale=1.5 ]  {$\boldsymbol \ldots$};
        

         \node[state,fill=pink] (An1) [right=25 mm of A_n] {$A_{n+1}'$};
    \node[state,fill=pink] (An11) [ below =25mm of An1]  {$A_{n+1}$};
    

\path
(A_1)  edge  [loop above=10]    node[]{\footnotesize{$a_0$}}     (A_1)
(A_1)  edge  node[xshift=0.1cm]{\footnotesize{$a_1$}}     (T_1)
(A_1)  edge  node[xshift=-0.2cm]{\footnotesize{$a_2$}}     (F_1)

(A_11)  edge  [loop above=10]    node[]{\footnotesize{$a_0$}}     (A_11)
(A_11)  edge  node[xshift=0.1cm]{\footnotesize{$a_1$}}     (T_11)
(A_11)  edge  node[xshift=-0.2cm]{\footnotesize{$a_2$}}     (F_11)

(T_1)  edge  node[xshift=-0.2cm]{\footnotesize{$a_3$}}     (A_2)
(T_1)  edge  node[xshift=-0.6cm,yshift=1.1cm]{\footnotesize{$a_4$}}     (A_22)
(F_1)  edge  node[xshift=-0.5cm,yshift=0.2cm,below]{\footnotesize{$a_3$}}     (A_2)
(F_1)  edge  node[xshift=-0.6cm, yshift=0.5cm]{\footnotesize{$a_4$}}     (A_22)
(T_11)  edge  node[below,xshift=-0.2cm]{\footnotesize{$a_3$}}     (A_22)
(F_11)  edge  node[below]{\footnotesize{$a_4$}}     (A_22)



(A_2)  edge[out=0, in=180, dashed]  node[xshift=-0.2cm,yshift=0.3cm,below]{\footnotesize{$a_5$}}     (T_2)
(A_2)  edge[out=0,in= 180,dashed]  node[below]{}     (F_2)

(A_22)  edge[out=0, in=180, dashed]  node[xshift=-0.2cm,yshift=0.3cm,below]{\footnotesize{$a_5$}}     (T_22)
(A_22)  edge[out=0,in= 180,dashed]  node[below]{}     (F_22)


(T_2)  edge  node[xshift=-0.2cm]{\footnotesize{$a_3$}}     (A_3)
(T_2)  edge  node[xshift=-0.6cm,yshift=1.1cm]{\footnotesize{$a_4$}}     (A_33)
(F_2)  edge  node[xshift=-0.5cm,yshift=0.2cm,below]{\footnotesize{$a_3$}}     (A_3)
(F_2)  edge  node[xshift=-0.6cm, yshift=0.5cm]{\footnotesize{$a_4$}}     (A_33)
(T_22)  edge  node[below,xshift=-0.2cm]{\footnotesize{$a_3$}}     (A_33)
(F_22)  edge  node[below]{\footnotesize{$a_4$}}     (A_33)



(A_3)  edge  [loop above=10]    node[]{\footnotesize{$a_0$}}     (A_3)
(A_3)  edge  node[xshift=0.1cm]{\footnotesize{$a_1$}}     (T_3)
(A_3)  edge  node[xshift=-0.2cm]{\footnotesize{$a_2$}}     (F_3)

(A_33)  edge  [loop above=10]    node[]{\footnotesize{$a_0$}}     (A_33)
(A_33)  edge  node[xshift=0.1cm]{\footnotesize{$a_1$}}     (T_33)
(A_33)  edge  node[xshift=-0.2cm]{\footnotesize{$a_2$}}     (F_33)

(T_3)  edge  node[xshift=-0.2cm]{\footnotesize{$a_3$}}     (A_4)
(T_3)  edge  node[xshift=-0.6cm,yshift=1.1cm]{\footnotesize{$a_4$}}     (A_44)
(F_3)  edge  node[xshift=-0.5cm,yshift=0.2cm,below]{\footnotesize{$a_3$}}     (A_4)
(F_3)  edge  node[xshift=-0.6cm, yshift=0.5cm]{\footnotesize{$a_4$}}     (A_44)
(T_33)  edge  node[below,xshift=-0.2cm]{\footnotesize{$a_3$}}     (A_44)
(F_33)  edge  node[below]{\footnotesize{$a_4$}}     (A_44)



(A_n)  edge[out=0, in=180, dashed]  node[xshift=-0.2cm,yshift=0.3cm,below]{\footnotesize{$a_5$}}     (T_n)
(A_n)  edge[out=0,in= 180,dashed]  node[below]{}     (F_n)

(A_n1)  edge[out=0, in=180, dashed]  node[xshift=-0.2cm,yshift=0.3cm,below]{\footnotesize{$a_5$}}     (T_n1)
(A_n1)  edge[out=0,in= 180,dashed]  node[below]{}     (F_n1)


(T_n)  edge  node[xshift=-0.2cm]{}     (An1)
(T_n)  edge  node[xshift=-0.7cm,yshift=1.1cm]{}     (An11)
(F_n)  edge  node[xshift=-0.4cm,yshift=-0.2cm,below]{}     (An1)
(F_n)  edge  node[xshift=-0.6cm, yshift=0.5cm,below]{}     (An11)
(T_n1)  edge  node[below,xshift=-0.2cm]{}     (An11)
(F_n1)  edge  node[below]{}     (An11)

(An1)  edge  [loop right=10]    node[]{\footnotesize{$a_0$}}     (An1)
(An11)  edge  [loop right=10]    node[]{\footnotesize{$a_0$}}     (An11)

;
\end{tikzpicture}}
\caption{An illustration of the MDP constructed to show that the BMP is PSPACE-hard. Solid lines represent deterministic transitions; dashed lines represent transitions with equal probability. The states $A_i',A_i$ represent the variables $x_i$ in the Boolean formula $F(x_1,x_2,\ldots,x_n)$. An existential variable $x_i$ is set to true (false) when the action $a_1$ ($a_2$) is incentivized from the states $A_i',A_i$.
A clause $C_k$ becomes true when the agent type $\theta_k$ transitions to a state $A_i$ with probability 1 under the incentive offers. The formula $F$ becomes true when the agent transitions to a state $A_i$ with probability 1 regardless of its type.  }
\label{fig:pspace_proof}
\end{figure*}
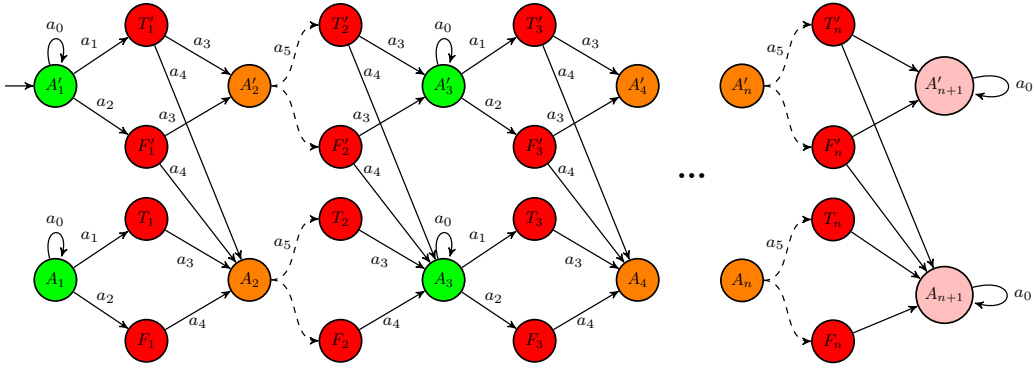

We first construct the MDP $\mathcal{M}$. A graphical illustration of the constructed MDP is given in Fig. \ref{fig:pspace_proof}. The set $S$ of states consists of $6n$$+$$2$ states. That is, $6$ states $A_i,A_i',T_i,T_i',F_i, F_i'$ for each variable $x_i$ where $i$$\in$$[n]$, and two additional states $A_{n+1},A_{n+1}'$. The initial state $s_0$ is the state $A_1'$. The set $\mathcal{A}$ of actions consists of 6 actions $a_0,a_1,\ldots,a_5$. The target set is $B$$=$$\{A_{n+1}\}$, and the type set is $\Theta$$=$$\{\theta_1,\theta_2,\ldots,\theta_m\}$.


\textit{Transition function $\mathcal{P}$:} We first define the transitions from the states $T_i,F_i,T_i',F_i'$. From $T_i$, under the only available action $a_3$, the agent makes a deterministic transition to $A_{i+1}$. From $F_i$, under the only available action $a_4$, the agent makes a deterministic transition to $A_{i+1}$. From $T_i'$, under the available actions $a_3$ and $a_4$, the agent deterministically transitions, respectively, to the states $A_{i+1}'$ and $A_{i+1}$. Similarly, from $F_i'$, under the available actions $a_3$ and $a_4$, the agent deterministically transitions, respectively, to the states $A_{i+1}'$ and $A_{i+1}$.

Next, we define the transitions from the states $A_i,A_i'$. From $A_i$, if $i$$\in$$[n]$ is even, the agent transitions to the states $T_i$ and $F_i$ with equal probability under the only available action $a_5$; if $i$$\in$$[n]$ is odd, the agent deterministically transitions to the states $A_i$, $T_i$, and $F_i$, under the actions $a_0$, $a_1$, and $a_2$, respectively. Similarly, from $A_i'$, if $i$$\in$$[n]$ is even, the agent transitions to the states $T_i'$ and $F_i'$ with equal probability under the only available action $a_5$; if $i$$\in$$[n]$ is odd, the agent deterministically transitions to the states $A_i'$, $T_i'$, and $F_i'$, under the actions $a_0$, $a_1$, and $a_2$, respectively. Finally, from the states $A_{n+1}$ and $A_{n+1}'$, the agent makes a self transition under the action $a_0$.

 \textit{Reward functions:} The reward function $\mathcal{R}_{\theta_k}$ where $k$$\in$$[m]$ is defined as follows. For $s$$\in$$\{A_i,A_i' :  i$$\in$$[n]\ \text{is odd}\}$, i.e., the green states in Fig. \ref{fig:pspace_proof}, we have
$\mathcal{R}_{\theta_k}(s,a_0)$$=$$0$, $\mathcal{R}_{\theta_k}(s,a_1)$$=$$-1$, $\mathcal{R}_{\theta_k}(s,a_2)$$=$$-1$. That is, the agent stays in the same state in the absence of incentives, and to move the agent toward the target state, the principal needs to offer the incentive of at least $1$$+$$\epsilon$ for $a_1$ or $a_2$, where $\epsilon$$>$$0$ is an arbitrarily small constant.


We now define the reward function for the states $T_i',F_i'$. This is the most important step in our construction since it relates the agent type $\theta_k$ with the clause $C_k$. If the variable $x_i$ appears positively in $C_k$ \footnote{A variable $x_i$ is said to appear positively in a disjunction clause $C_k$ if $C_k$ is true when $x_i$ is true; it is said to appear negatively in $C_k$ if $C_k$ is true when $x_i$ is false; it is said to not appear in $C_k$ if the truth value of $C_k$ is independent of the truth value of $x_i$. For example, in $C$$=$$x_1 \lor \lnot x_2$, the variable $x_1$ appears positively, $x_2$ appears negatively, and $x_3$ does not appear.}, for $s$$\in$$\{T_i': i$$\in$$[n]\}$, we have 
\begin{align*}
    \mathcal{R}_{\theta_k}(s,a_3)=n+1, \ \mathcal{R}_{\theta_k}(s,a_4)=0,
\end{align*}
and for $s$$\in$$\{F_i': i$$\in$$[n]\}$, we have 
\begin{align*}
    \mathcal{R}_{\theta_k}(s,a_3)=0, \ \mathcal{R}_{\theta_k}(s,a_4)=n+1.
\end{align*}

On the other hand, if the variable $x_i$ appears negatively in $C_k$, for $s$$\in$$\{T_i': i$$\in$$[n]\}$, we have 
\begin{align*}
    \mathcal{R}_{\theta_k}(s,a_3)=0, \ \mathcal{R}_{\theta_k}(s,a_4)=n+1,
\end{align*}
and for $s$$\in$$\{F_i': i$$\in$$[n]\}$, we have 
\begin{align*}
    \mathcal{R}_{\theta_k}(s,a_3)=n+1, \ \mathcal{R}_{\theta_k}(s,a_4)=0.
\end{align*}
Finally, if the variable $x_i$ does not appear in $C_k$, for $s$$\in$$\{T_i': i$$\in$$[n]\}$, we have
\begin{align*}
    \mathcal{R}_{\theta_k}(s,a_3)=0, \ \mathcal{R}_{\theta_k}(s,a_4)=n+1,
\end{align*}
and for $s$$\in$$\{F_i': i$$\in$$[n]\}$, we have 
\begin{align*}
    \mathcal{R}_{\theta_k}(s,a_3)=0, \ \mathcal{R}_{\theta_k}(s,a_4)=n+1.
\end{align*}
For all the other states and actions, we have $\mathcal{R}_{\theta_k}(s,a)=0$.


We now show that an incentive sequence that leads the agent to the target state $A_{n+1}$ with probability 1 at the worst-case total cost $n$ or less exists if and only if the QBF $F$ is true.

Suppose first that the QBF $F$ is true. Then, there is a truth value assignment, i.e., true or false, for each existential variable $\{x_i: i\in[n] \ \text{is odd}\}$, such that under this assignment all clauses $C_k$ where $k$$\in$$[m]$ are true. Then, on the constructed MDP, from the states $\{A_i': i\in[n] \ \text{is odd}\}$, i.e., the green states in Fig. \ref{fig:pspace_proof}, we can simply offer the following stationary incentives to the agent:
If $x_i$ is true, $\gamma(A_i',a_1)$$=$$2$; if $x_i$ is false $\gamma(A_i',a_2)$$=$$2$. For all the other state-action pairs $(t,a)$, we offer $\gamma(t,a)$$=$$0$. Since the formula $F$ is true, all the clauses $C_k$ must be true. Then, under the provided incentives, each type $\theta_{k}$ must eventually transition to a state $A_i$ with probability 1. Since, from any state $A_i$, the agent reaches to the target state $A_{n+1}$ with probability 1, the reachability constraint is satisfied. Moreover, the principal pays the incentives $\gamma(A_i',a_1)$$=$$2$ exactly $n/2$ times; hence, the total cost to the principal is $n$. 

Suppose that there exists an incentive sequence under which the agent reaches the state $A_{n+1}$ with probability 1 at the worst-case total cost $n$ or less. Then, the optimal policy of the agent under such an incentive sequence differs from its optimal policy in the absence of incentives only in the states $\{A_i': i$$\in$$[n] \ \text{is odd}\}$. The previous claim is true because changing the optimal policy of the agent in states $T_i',F_i'$ requires the principal to offer at least an incentive amount equal to $n+1$. Note that all the agent types $\theta_k$ reaches the target state with probability 1 under the provided incentive sequence. Then, it must be true that the principal incentivizes the actions $a_1$ and/or $a_2$ in states  $\{A_i': i\in[n] \ \text{is odd}\}$ such that the agent eventually reaches a state $A_i$ with probability 1 regardless of its type. We set the existential variable $x_i$ to true if the principal incentives $a_1$, and to false if the principal incentivizes only $a_2$. Recall that each clause $C_k$ corresponds to an agent type $\theta_k$ and the transition of an agent type $\theta_k$ to a state $A_i$ with probability 1 implies the clause $C_k$ becoming true.  Consequently, since the agent reaches the state $A_{n+1}$ with probability 1 \textit{regardless of its type}, under the described truth value assignment, the formula $F$ is evaluated true. $\Box$

{\setlength{\parindent}{0cm}\begin{definition} (Euclidian path-TSP) \cite{papadimitriou1977euclidean} Given a set $[N]$ of cities, distances $c_{i,j}$$\in$$\mathbb{N}$ between each city pair $(i,j)$ such that $c_{i,j}$$=$$c_{j,i}$ and $c_{i,j}$$+$$c_{j,k}$$\geq$$c_{i,k}$ for all $i,j,k$$\in$$[N]$, and a constant $K$$\in$$\mathbb{N}$, decide whether there exists a path from the city 1 to city N that visits all cities exactly once and the total traversed distance is $K$ or less.
\end{definition}}

\noindent \textbf{Proof of Theorem \ref{N-BMP_thm}:} The decision problem is in NP because, for a given stationary incentive sequence, we can compute the occupancy measure \cite{Altman} for the Markov chain induced by the optimal \textit{stationary} policy of each agent type in polynomial-time via matrix inversion. The reachability probability and the incurred total cost are linear functions of the occupancy measure \cite{Puterman}. Hence, the satisfaction of the reachability constraint and the corresponding total cost incurred by the principal can be verified in polynomial-time. 

The NP-hardness proof is by a reduction from the Euclidian path-TSP problem. We are given an arbitrary Euclidian path-TSP problem instance. We first construct an MDP $\mathcal{M}$. The set of states is $S$$=$$\{q_1,q_2,\ldots,q_N\}$, the initial state is $q_1$, and the set of actions is $\mathcal{A}$$=$$\{a_1,a_2,\ldots,a_N\}$. The transition function $\mathcal{P}$ is such that, for states $q_i$$\in$$S\backslash \{q_N\}$, $\mathcal{P}_{q_i,a_j,q_j}$$=$$1$ for all $j$$\in$$\{1,2,\ldots,N\}$ and, the state $q_N$$\in$$S$ is absorbing. 

The type set is $\Theta$$=$$\{\theta_1,\theta_2,\ldots,\theta_{N-1}\}$, and the target set is $B$$=$$\{q_N\}$. The reward function $\mathcal{R}_{\theta_i}$ for the type $\theta_i$$\in$$\Theta$ is
\begin{align*}
    \mathcal{R}_{\theta_i}(q_j,a_k)=
    \begin{cases} 
    0 & \text{if}\  j=k\\
    -c_{j,k} & \text{if} \  j\neq k,\  j=i\\
    -c_{j,k} & \text{if} \  j\neq k, \ j\neq i,\ k\neq N \\
    -(K+1) & \text{if} \  j\neq k,\ j\neq i, k=N.
    \end{cases}
\end{align*}

We claim that there exists a feasible solution to the constructed N-BMP with the objective value $K$ or less if and only if there exists a TSP path with total cost $K$ or less. 


Suppose that there exists a feasible solution to the N-BMP with the objective value of $K$ or less. Then, there exists an incentive sequence under which at least one agent type visits all states in $S$ exactly once. Suppose for contradiction that none of the agent types visits the state $q_i$$\in$$S$ and that the true agent type is $\theta^{\star}$$=$$\theta_i$. Then, to convince the agent to reach the target state $q_N$, the principal must pay the incentive of $(K+1)$ for the action $a_N$ from at least one of the states $q$$\in$$S\backslash \{q_i\}$, because otherwise, the agent will never take the action $a_N$ and reach the target state. Therefore, the optimal incentive sequence incurs the total cost of at least $(K$$+$$1)$, which raises a contradiction. Moreover, the agent type $\theta_i$ cannot visit the same state twice under a provided \textit{stationary} incentive sequence because otherwise the agent's \textit{stationary} optimal policy violates the reachability constraint. Consequently, the path followed by the agent type $\theta_i$ constitutes a solution to the Euclidian path-TSP with the total cost $K$ or less.

Consider a solution to the Euclidian path-TSP problem which visits the states $s_1,s_2,\ldots,\text{and}\ s_N$ in sequence, where $s_1$$=$$q_1$ and $s_N$$=$$q_N$. The principal can construct a stationary incentive sequence from this solution as follows. For all states $s_i$ where $i$$\in$$[N]$, offer the incentives $\gamma(s_i,a_N)$$=$$c_{s_i,s_N}$$+$$\epsilon$ and $\gamma(s_i,a_{s_{i+1}})$$=$$c_{s_i,s_{i+1}}$$+$$\epsilon$. Under the constructed incentive sequence, each agent type $\theta_i$$\in$$\Theta$ transitions to the target state $q_N$ from the state $q_i$. Moreover, the worst-case total cost to the principal is equal to the cost of the Euclidian path-TSP problem plus $N\epsilon$, which can be made arbitrarily small. Hence, the claim follows. $\Box$
{\setlength{\parindent}{0cm}\begin{definition} (Set cover) \cite{karp1972reducibility}  Given a set $\mathcal{S}$, a collection of subsets $\mathcal{T}$$=$$\{T_i\subseteq \mathcal{S}\}$ such that $\cup_{T_i\in \mathcal{T}}$$=$$\mathcal{S}$, and a positive integer $M$, decide whether there exists a subcollection $\mathcal{U}$$\subseteq$$\mathcal{T}$ of size $M$ or less, i.e., $\lvert \mathcal{U} \rvert$$\leq$$M$, that covers $\mathcal{S}$, i.e., $\cup_{T_i\in \mathcal{U}}$$=$$\mathcal{S}$.
\end{definition}}

\noindent \textbf{Proof of Theorem \ref{NS-BMP_thms}:} The decision problem can be shown to be in NP using the same method presented in Theorem \ref{N-BMP_thm}. The NP-hardness proof is by a reduction from the set cover problem. We are given an arbitrary set cover problem instance in which $\mathcal{S}$$=$$\{1,2,\ldots,N\}$ and $\mathcal{T} $$=$$\{T_1,T_2,\ldots,T_K\}$ such that $M$$\leq$$K$. To prove the claim, we construct an MDP $\mathcal{M}$, a target set $B$$\subseteq$$S$, a type set $\Theta$, and a state-independent reward function $\mathcal{R}_{\theta}$ for each $\theta$$\in$$\Theta$. Then, we show that, on the constructed model, the total cost of an incentive sequence that solves the NS-BMP is $M$ or less if and only if $\rvert \mathcal{U}\lvert$$\leq$$M$.

The MDP $\mathcal{M}$ has $M$$+$$2$ states, i.e., $S$$=$$\{q_1,q_2,\ldots,q_{M+2}\}$, and $K$$+$$1$ actions, i.e., $\mathcal{A}$$=$$\{a_0,a_1,\ldots,a_{K}\}$. The initial state is $q_1$. From the states $q_i$ where $i$$\in$$[M$$-$$1]$, the agent transitions to the state $q_{i+1}$ under the action $a_0$, and to the state $q_{M+1}$ under all the other actions. From the state $q_{M}$, the agent transitions to the state $q_{M+2}$ under the action $a_0$, and to $q_{M+1}$ under all the other actions. The states $q_{M+1}$ and $q_{M+2}$ are absorbing.

We define $B$$=$$\{q_{M+1}\}$ and  $\Theta$$=$$\{\theta_1,\theta_2,\ldots,\theta_N\}$. Finally, the state-independent reward function $\mathcal{R}_{\theta_i}$ is defined as
\begin{align*}
    \mathcal{R}_{\theta_i}(q,a_j)=\begin{cases}0 & \text{if}\ j=0\\
    -1/2 & \text{if} \ i\in T_j\\
    -(K+1) & \text{otherwise}.
    \end{cases}
\end{align*}
Note that, in the absence of incentives, the agent reaches the absorbing state $q_{M+2}$ with probability 1 regardless of its type.

Suppose there exists a collection of subsets $\mathcal{U}$$\subseteq$$\mathcal{T}$ of size $M$ or less such that $\cup_{T_i\in \mathcal{U}}$$=$$\mathcal{S}$. Without loss of generality, assume that $\mathcal{U}$$=$$\{T_1,T_2,\ldots,T_L\}$ where $L$$\leq$$M$. Consider the following stationary deterministic incentive sequence
\begin{align*}
    \gamma(q_i,a_j)=\begin{cases} 1 & \text{if} \ i=j\\
    0 &\text{otherwise}.
    \end{cases}
\end{align*}
Under the incentive sequence given above, an agent type $\theta_i$ transitions to the state $q_{M+1}$ from a state $q_k$ where $k$$\in$$[M]$ if $i$$\in$$\cup_{T_j\in\mathcal{U}}$. Since $\cup_{T_j\in\mathcal{U}}$$=$$\mathcal{S}$, all agent types transitions to the target state $q_{M+1}$ with probability 1. Moreover, the total cost to the principal is clearly less than or equal to $M$. 

Suppose there exists a stationary deterministic incentive sequence that leads the agent to the target state $q_{M+1}$ regardless of its type at the total cost of $M$ or less. Then, the incentive sequence is such that, the agent type $\theta_i$ where $i$$\in$$T_j$ transitions to the target state $q_{M+1}$ from a state $q_k$ where $k$$\in$$[M]$ if the action $a_j$ is incentivized from that state. Since the agent reaches the target state regardless of its type, the collection of subsets $T_j$ that corresponds to the incentivized actions $a_j$ constitutes a set cover. Because the incentives are provided only from the states $q_k$ where $k$$\in$$[M]$, the size of the resulting set cover is less than or equal to $M$. $\Box$

\noindent\textbf{Proof of Theorem \ref{theorem_dominant}:} For a given incentive sequence $\gamma$$\in$$\Gamma_{R,\theta_d}(\mathcal{M})$, let $\pi^{\star}$$=$$(d_1^{\star},d_2^{\star},\ldots)$ be the optimal policy of the agent type $\theta_d$ and $\overline{\gamma}$$\in$$\Gamma(\mathcal{M})$ be an incentive sequence
\begin{align}\label{dom_proof_1}
    \overline{\gamma}(I_t,a)=\begin{cases} \gamma(I_t,a) & \text{if} \ d_t^{\star}(s)=a\\
    0 &\text{otherwise}.
    \end{cases}
\end{align}
Note that $\mathcal{R}_{\theta_d}(s,a)$$+$$\overline{\gamma}(I_t,a)\geq\max_{a'\in \mathcal{A}(s)} \mathcal{R}_{\theta_d}(s,a')$ for $d_t^{\star}(s)$$=$$a$ and $\overline{\gamma}$$\in$$\Gamma_{R,\theta_d}(\mathcal{M})$. The condition stated in the theorem implies that
    $\mathcal{R}_{\theta}(s,a)$$+$$\overline{\gamma}(I_t,a)\geq\max_{a'\in \mathcal{A}(s)} \mathcal{R}_{\theta}(s,a')$
for any $\theta$$\in$$\Theta$. Hence, we have $\overline{\gamma}$$\in$$\Gamma_{R,\Theta}(\mathcal{M})$, which implies that $\Gamma_{R,\Theta}(\mathcal{M})$$\subseteq$$\Gamma_{R,\theta_d}(\mathcal{M})$. Since, $\Gamma_{R,\theta_d}(\mathcal{M})$$\subseteq$$\Gamma_{R,\Theta}(\mathcal{M})$ by definition, we conclude that $\Gamma_{R,\Theta}(\mathcal{M})$$=$$\Gamma_{R,\theta_d}(\mathcal{M})$.

Now, for any given incentive sequence $\gamma$$\in$$\Gamma_{R,\Theta}(\mathcal{M})$, let $\pi^{\star}$$=$$(d_1^{\star},d_2^{\star},\ldots)$ be the optimal policy of the agent type $\theta_d$ and $\overline{\gamma}$$\in$$\Gamma(\mathcal{M})$ be an incentive sequence defined as in \eqref{dom_proof_1}. If the condition stated in the theorem holds, then we have 
\begin{align*}
    \max_{\theta\in\Theta}f(\overline{\gamma},\theta)= f(\overline{\gamma},\theta_d).
\end{align*}
Using the identity $\Gamma_{R,\Theta}(\mathcal{M})$$=$$\Gamma_{R,\theta_d}(\mathcal{M})$, we conclude that 
\begin{align*}
    \min_{\gamma\in\Gamma_{R,\Theta}(\mathcal{M})}\max_{\theta\in\Theta}f(\gamma,\theta)= \min_{\gamma\in\Gamma_{R,\theta_d}(\mathcal{M})}f(\gamma,\theta_d).\  \Box
\end{align*}

\section{}\label{appendix_algorithm}
In this appendix, we show how to modify the algorithms presented in Section \ref{algorithms_section} in order to compute globally and locally optimal solutions to the NS-BMP.

Let $\{\overline{X}_{s,a}$$\in$$\{0,1\}: s$$\in$$S_r, a$$\in$$\mathcal{A}, \sum_{a\in\mathcal{A}}\overline{X}_{s,a}$$=$$1 \}$ be a set of binary variables. One can obtain a globally optimal solution to the NS-BMP by adding the constraints 
\begin{align}\label{additional_cons_app}
    \sum_{a\in\mathcal{A}}\overline{X}_{s,a}\gamma(s,a)=\sum_{a\in\mathcal{A}}\gamma(s,a)\ \  \text{ for all $s$$\in$$S_r$}
\end{align}
to the MILP \eqref{MILP_begin}-\eqref{MILP_end} and solving the resulting optimization problem. Note that the above constraint is not linear in the variables $\gamma(s,a)$; however, each term $\overline{X}_{s,a}\gamma(s,a)$ can be replaced exactly by its corresponding McCormick envelope in order to obtain an MILP formulation for the NS-BMP. 

Let $\{\overline{\nu}_{s,a}$$\geq$$0: s$$\in$$S_r, a$$\in$$\mathcal{A}, \sum_{a\in\mathcal{A}}\overline{\nu}_{s,a}$$=$$1 \}$ be a set of continous variables. One can obtain a locally optimal solution to the NS-BMP by adding the constraints 
\begin{align}
    \sum_{a\in\mathcal{A}}\overline{\nu}_{s,a}\gamma(s,a)\geq\sum_{a\in\mathcal{A}}\gamma(s,a)\ \  \text{ for all $s$$\in$$S_r$}
\end{align}
to the NLP formulated in Section \ref{local_opt_sec} and solving the resulting optimization problem by resorting to the CCP.

\vspace{-1.2cm}
\begin{IEEEbiography}[{\includegraphics[width=1in,height=1.25in,clip,keepaspectratio]{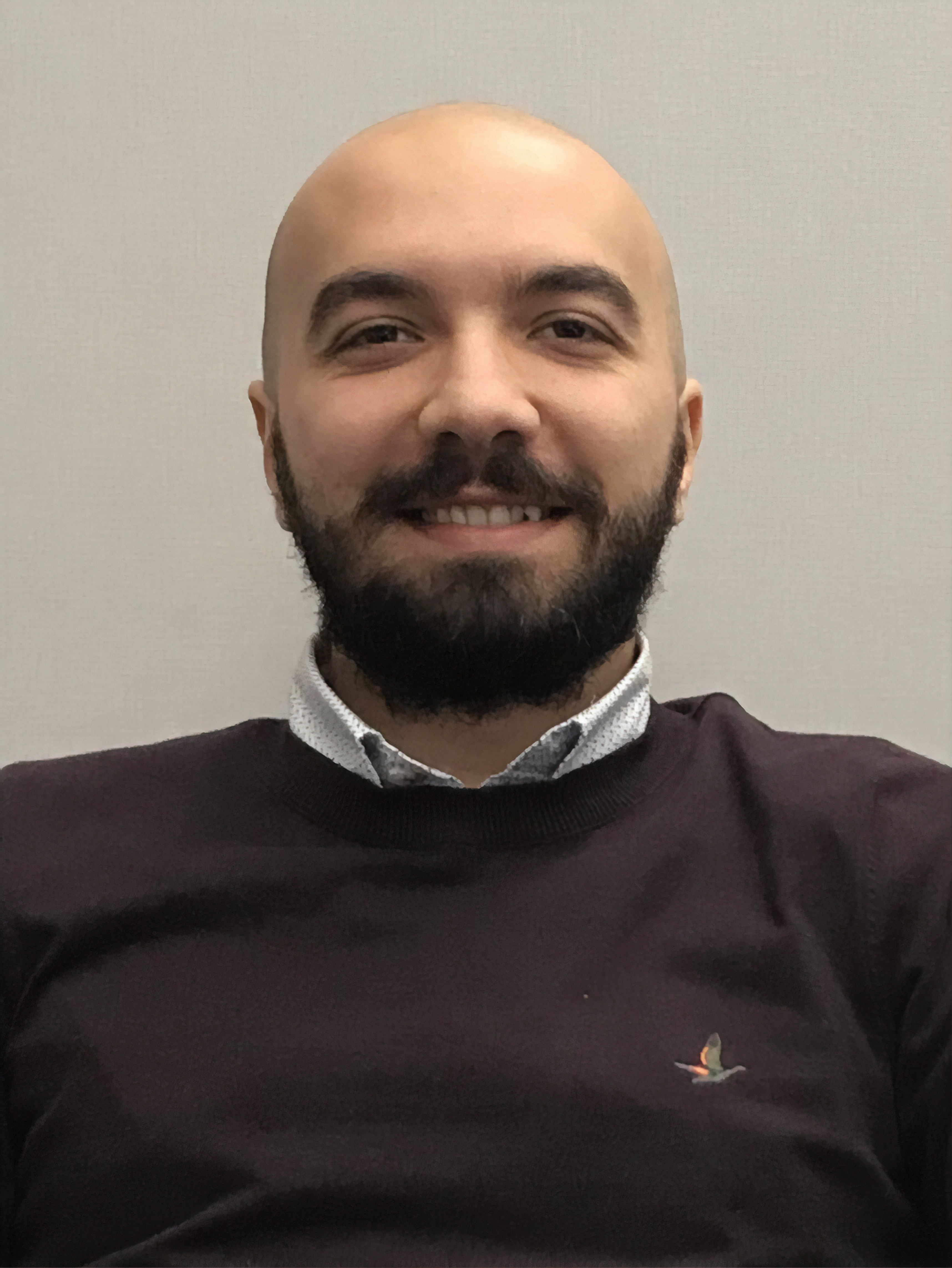}}]{Yagiz Savas} joined the Department of Aerospace Engineering at the University of Texas at Austin as a Ph.D. student in Fall 2017. He received his B.S. degree in Mechanical Engineering from Bogazici University in 2017. His research focuses on developing theory and algorithms that guarantee desirable behavior of autonomous systems operating in uncertain, adversarial environments.
\end{IEEEbiography}\vspace{-1.2cm}
\begin{IEEEbiography}[{\includegraphics[width=1in,height=1.25in,clip,keepaspectratio]{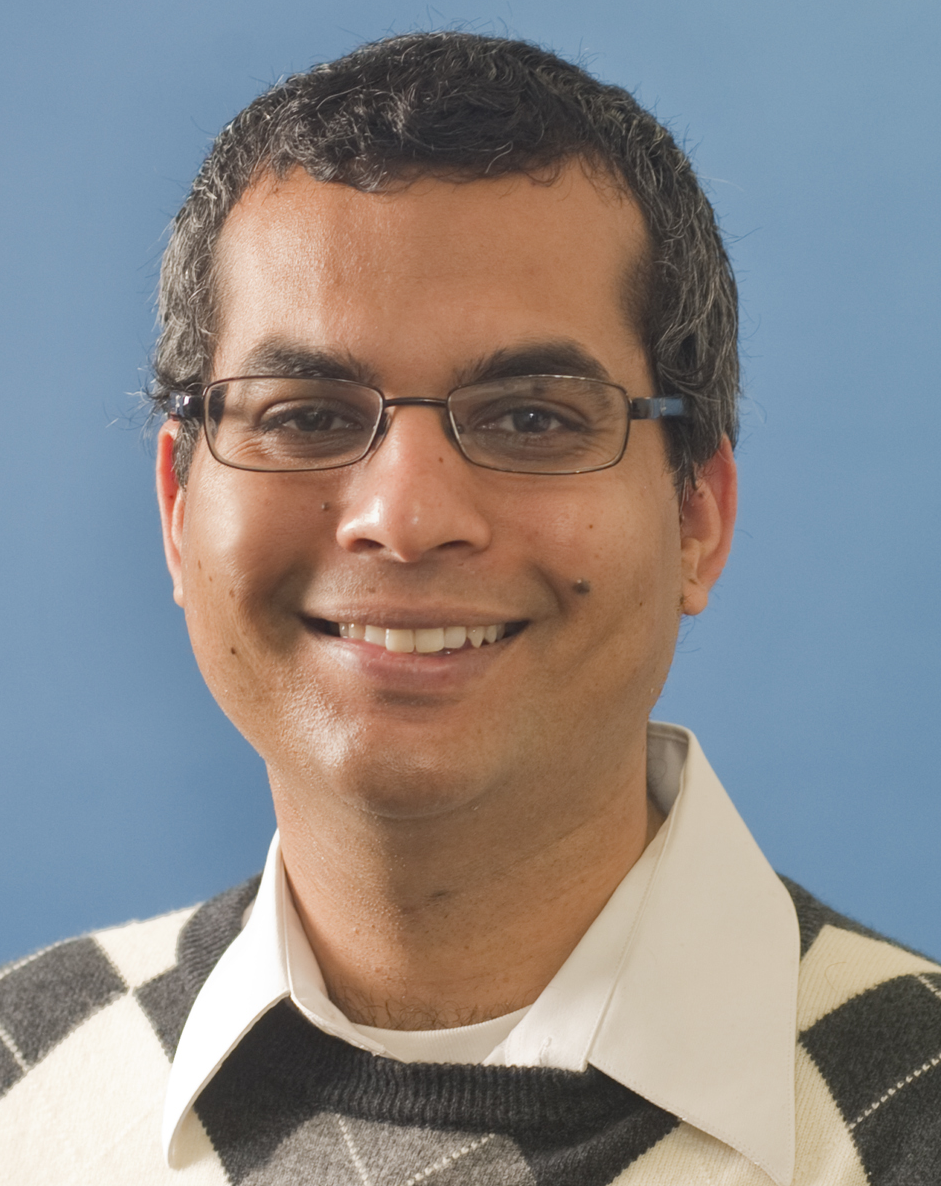}}] {Vijay Gupta} is in the Department of Electrical Engineering at the University of Notre Dame since 2008.  He received the 2018 Antonio J Rubert Award from the IEEE Control Systems Society, the 2013 Donald P. Eckman Award from the American Automatic Control Council and a 2009 National Science Foundation (NSF) CAREER Award. His research interests are broadly in the interface of communication, control, distributed computation, and human decision making.
\end{IEEEbiography}\vspace{-1.2cm}
\begin{IEEEbiography}[{\includegraphics[width=1in,height=1.25in,clip,keepaspectratio]{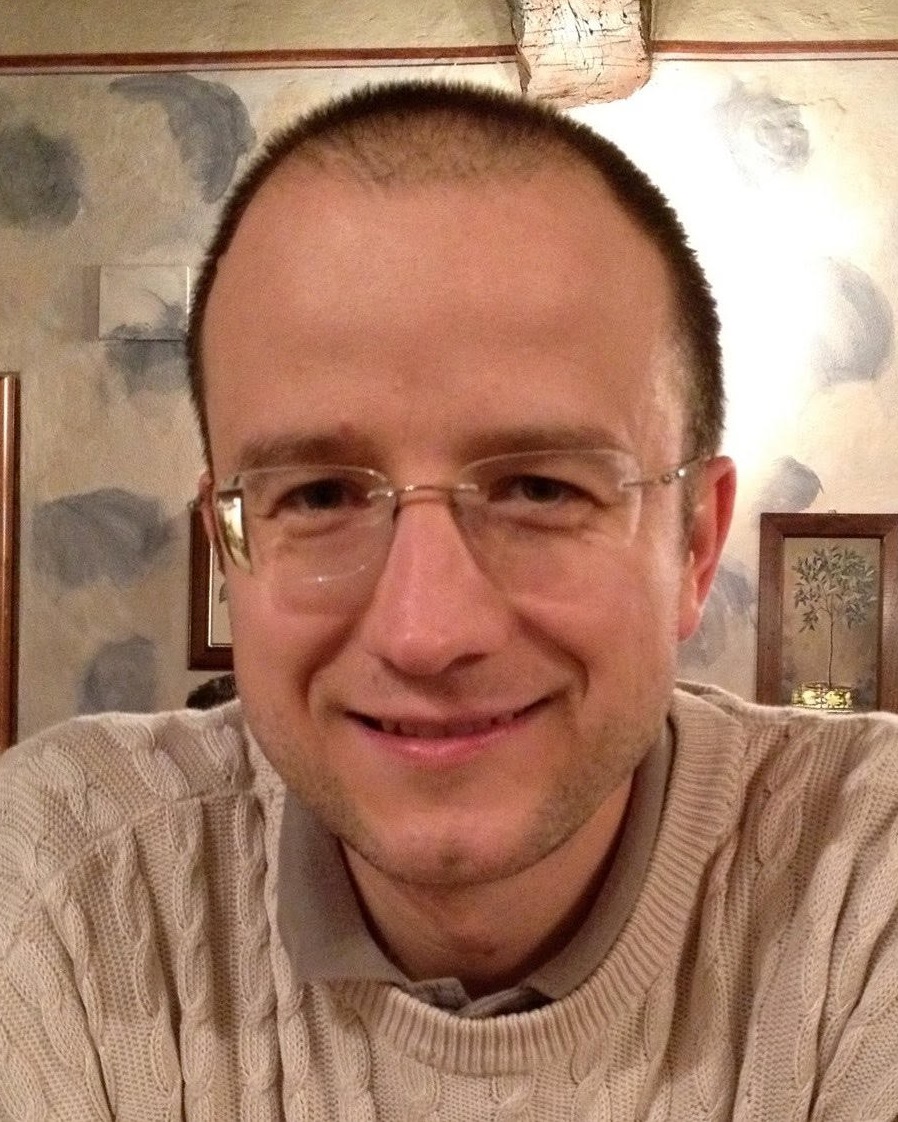}}]{Ufuk Topcu} joined the Department of Aerospace Engineering at the University of Texas at Austin as an assistant professor in Fall 2015. He received his Ph.D. degree from the University of California at Berkeley in 2008. He held research positions at the University of Pennsylvania and California Institute of Technology. His research focuses on the theoretical, algorithmic and computational aspects of design and verification of autonomous systems through novel connections between formal methods, learning theory and controls.
\end{IEEEbiography}
\end{document}